\documentclass[10pt]{article}

\usepackage{style}
\newcommand{\aref}[1]{Appendix~\ref{#1}}
\newcommand{\sref}[1]{Section~\ref{#1}}
\newcommand{\fref}[1]{Figure~\ref{#1}}

\title{\textbf{Likelihood estimation of an interpretable early-warning sign of critical transitions}}
\author[1]{Davide Papapicco}
\author[1]{Lauren Smith}
\author[1]{Graham Donovan}

\affil[1]{Department of Mathematics, University of Auckland, New Zealand}
\date{}

\begin{document}

\maketitle

\begin{abstract}
        A critical transition is a sudden, abrupt and potentially unforeseen change in the state of a complex system often associated with catastrophic collapses and irreversability.
        Here we introduce a novel, interpretable early-warning signal of bifurcation-induced critical transitions.
        The proposed indicator can be readily computed from timeseries observations by maximum likelihood estimation and it has the advantage of providing immediate insight on the state of the system.
        We showcase its applicability by assessing its performance using models of ecosystems and ocean circulation collapse.
\end{abstract}

\section{Introduction}
Over the past two decades much attention has been gathered in the study of catastrophic collapses of complex systems under time-varying forcing.
Initially described as \textit{``loss of resilience''} in ecosystems \cite{Ives1995, Scheffer2009}, climate models \cite{Lenton2012} and even social dynamics \cite{Scheffer2021}, these critical transitions, and their prediction, have been the subject of investigation across many disciplines \cite{Feng2014,Liang2017,Xu2023}.
These events, sometimes refered to as tipping points, are abrupt and sudden changes in the state of a system subject to external forcing.
More specifically, we say that a system has reached a tipping point when the observed steady-state evolves, over time, in a relatively stable fashion before suddenly \textit{``jumping''} to a different regime.

Different types of tipping events lead to qualitatively different behaviour of the system undergoing the transition, depending on the underlying mechanism that generates them. 
A classification of those was provided in \cite{Ashwin2012}:
\begin{itemize}
        \item bifurcation-induced (B-)tipping events are abrupt transitions of the steady-state once a critical threshold of a slowly changing forcing parameter is crossed. These types of tipping events are generically associated with saddle-node bifurcations;
     \item rate-induced (R-)tipping events are caused by the failure of the steady-state to track the drift of a quasi-steady equilibrium. The mechanism generating these tipping phenomena is fundamentally different from B-tipping as it is the rate of change at which the forcing parameter is ramped that causes the tipping \cite{Ashwin2017}, rather than its value. Notice that the steady-state need not undergo a global nor local bifurcation for R-tipping to occur;
     \item noise-induced (N-)tipping events are entirely due to stochastically driven jumps of the state of the system outside the basin attraction of the steady-state and are thus hard to predict. Notice that N-tipping do not need time-varying forcing of the parameter and can be observed in systems at full stationarity.
\end{itemize}

Here we focus on B-tipping events as it is believed they are the mechanics behind the catastrophic collapses in ecosystems \cite{May1977}, biology \cite{Donovan2015} and ocean currents \cite{Lohmann2024,vanWesten025a}.
A common characteristic of these systems is that they all exhibit multistability in which alternative stable equilibria coexist simultaneously.
The loss of stability or the vanishing of one such equilibria are modeled using bifurcation normal forms.
In a B-tipping event two equilibria of the system collide, annihilating each other.
A system tracking one of those equilibria will thus be forced to tip abruptly onto alternative stable states.
The new equilibrium state that follows the B-tipping can have catastrophic repercussions for the health of the system.

Due to the disruptive nature of tipping events, strategies have been developed that can predict these events by monitoring the state evolution.
Early-warning signals (EWS) are (often statistical) quantities that can be extracted from timeseries data with the intention of informing how likely and how close a critical transition is.
These precursors often rely on mechanistic features of the system approaching B-tipping of which \textit{critical slowing down} (CSD) has been the most widely studied \cite{VanNes2007, Dakos2008}.

CSD quantifies slower recovery to equilibrium from perturbations as the leading eigenvalue of the linearized system approaches the imaginary axis.
This property is remarkably generic as it has been shown that CSD can be detected in models of ecosystem collapse \cite{Scheffer2009,Dakos2012}, ocean circulation \cite{Boers2021a,Boers2021b}, climate and weather extreme events \cite{Lenton2012,Prettyman2019}, biological and physiological systems \cite{Donovan2015}, traffic flow \cite{Chattopadhyay2026} and many others.
Among the many proposed indicators, one of the most widely used methods to capture CSD is by detecting the increase in timeseries variance of the sample paths as the system is driven closer to B-tipping.

Despite the many efforts, EWS are still largely plagued by a general lack of interpretability \cite{Ditlevsen2010, Donovan2026}, unknown generalization to high-dimensional systems \cite{Robinson2025} and relative difficulty in assessing their robustness to real-world data \cite{Dakos2012, Ashwin2025}.
Let us consider the issue of interpretability using the increase in timeseries variance as an EWS (noting that the same issue arises for other EWS in the literature).
In an idealized scenario, the increase in variance is monotonic as the system approaches B-tipping and captures CSD. The question remains: how much increase in the variance is indicative of tipping?

Another way of framing this issue is the following: suppose we are given a single sample describing the time evolution of a given system of interest.
If we compute the variance of such sample, how is this single datum able to inform us on how close the system is to a tipping event?
To answer this question we would need to compare the computed variance to a baseline value for the same system when at a safe distance from the B-tipping event.
Moreover, different systems have different baseline values which makes the sample variance even more limiting as an EWS.
This limitation is especially significant if one only has access to limited data and no ensemble of timeseries is available, which is usually the case in the monitoring of real-world systems.

While in some case-specific applications interpretable EWS can be obtained, as demonstrated by studies for animal population collapse \cite{Drake2010}, disease emergence and resurgence \cite{Brett2020,Harris2020} and large scale ocean circulations \cite{Rodal2022}, a lack of universally interpretable, system-agnostic EWS still persists.
Specifically, the interpretability provided by these decision rules rely either on existing statistical indicators of proximity to B-tipping \cite{Drake2010,Harris2020} or on training datasets that can only generalize to similar systems \cite{Brett2020}.

In this paper we introduce a novel EWS that solves this issue. 
We achieve this by reframing the problem of detecting critical transitions in timeseries through the lens of statistical mechanics, leveraging known results from large deviations theory such as Kramer's formula.
Under mild assumptions, the proposed EWS is bounded in $(0,1]$ and monotonically increasing as the B-tipping is approached, regardless of the details of the system being monitored. 
Values of the proposed quantity that are close to $0$ translate to a system being far from a critical transition of B-tipping type, while values close to $1$ are indicative of proximity of a tipping event.
A full computational pipeline also accompanies the proposed new measure, which can be readily applied to input timeseries without prior knowledge of the model generating it.

This manuscript is organised in the following structure: in \sref{sec:escape_rate} we introduce the new EWS as a modified version of Kramer's law for the escape rate of overdamped particles across a potential barrier; 
in \sref{sec:mle} we discuss numerical strategies to extract the proposed EWS from timeseries data based on maximum likelihood estimators; 
this is followed in \sref{sec:applications} by comparisons of the performance of the proposed EWS with respect to the canonical increase in variance.
Applications in this section include a well-known model of vegetation biomass collapse in ecosystems, when subject to increasing harvesting, and an equally popular model of the Atlantic Meridional Overturning Circulation collapse under increasing freshwater forcing in the North Atlantic. 
In \sref{sec:analysis} we analyse the accuracy, sensitivity and statistical significance of the proposed estimator under relaxations of its working assumptions.
We conclude the paper with further discussions and future directions based on the presented technique in \sref{sec:discussion}.
\section{Escape rate across a potential barrier}\label{sec:escape_rate}
Critical transitions in complex systems have long been modeled as qualitative changes in the solutions of multiscale stochastic differential equations (SDEs) subject to slow forcing \cite{Kuehn2011}.
For finite-dimensional systems, these are usually presented in the form
\begin{equation}
        \begin{aligned}\label{eq:fast_slow_sde}
                dX_{t} &= f(X_{t},\mu)dt + \sigma dW_{t}\,,\\
                d\mu &= \varepsilon dt\,,
        \end{aligned}
\end{equation}
where $X_t$ denotes a random variable (r.v.) modeling the state of the system of interest in $\mathbb{R}^{d}$ while $\mu\in \mathbb{R}^{p}$ denotes a parameter vector.
The time evolution of the r.v. follows a deterministic drift, described by the parametrized vector field $f:\mathbb{R}^{d}\times \mathbb{R}^{p}\to \mathbb{R}^{d}$, and a stochastic diffusion of intensity $\sigma\in\mathbb{R}$, modeled by a $d-$dimensional Wiener process $W_{t}$ with zero mean and Gaussian distributed increments (i.e. $W_{s} - W_{t}\sim\mathcal{N}(0,s-t)$).
For simplicity, we will assume that the stochastic fluctuations are state-independent (i.e. \eqref{eq:fast_slow_sde} is an additive noise SDE in $X_t$).
Importantly, the parameter characterizing the vector field of the dynamics is not stationary but it evolves according to a slow (deterministic) linear ramp of slope $0 < \varepsilon \ll 1$.
This setup allows us to slowly change the state of the system, hence reflecting the multiscale time-varying forcing of real-world complex systems.

In the following, and for the rest of this paper, we will treat functions $y$ that depend both on the state variable $x$ and on time $t$, i.e. $y=y(x,t)$.
To avoid confusion with regards to differentiation we will adopt the following notation
\begin{equation*}
        \dot{y}:=\partial_t\,y(x,t)\quad \text{and}\quad y':=\partial_x\,y(x,t)\,.
\end{equation*}
We will also assume that, for fixed parameters $\mu\in \mathbb{R}^{p}$ the (frozen) deterministic drift is a conservative vector field, i.e. there exists a scalar potential $V(x;\mu):\mathbb{R}^{d}\to \mathbb{R}$ such that $f(x;\mu) = -V'(x;\mu)$.
Since we would like to predict B-tipping of saddle-node (SN) type, we assume that the multistability of the system entails that $V$ is a sufficiently smooth function of the state variable $x$ and that it also grows sufficiently fast as $x\to\pm\infty$.
The gradient structure of the dynamics allows to map equilibria of the vector field $f$ to stationary points of its potential $V$.
In particular, locally stable equilibria of $f$ correspond to local minima of $V$ while locally unstable equilibria are associated to local maxima of $V$.
This allows for a useful interpretation of the solution of \eqref{eq:fast_slow_sde} as the motion of an overdamped particle governed by the Langevin equation
\begin{equation}\label{eq:langevin}
        dX_{t} = -V'(X_{t};\mu)dt + \sqrt{2D}dW_{t}\,,\quad D=\frac{\sigma^{2}}{2}\,.
\end{equation}
If the process is stationary (i.e. $\varepsilon=0$) then, for a given configuration of the potential $V$ at parameter value $\mu$, the long term behaviour of a solution of the Langevin equation of motion \eqref{eq:langevin} is described by random fluctuations (of intensity $\sqrt{2D}$) around the local minima of $V$.
\subsection{Kramer's escape law}\label{subsec:kramer_escape_formula}
\begin{figure}[!t]
	      \includegraphics[keepaspectratio, width=\textwidth]{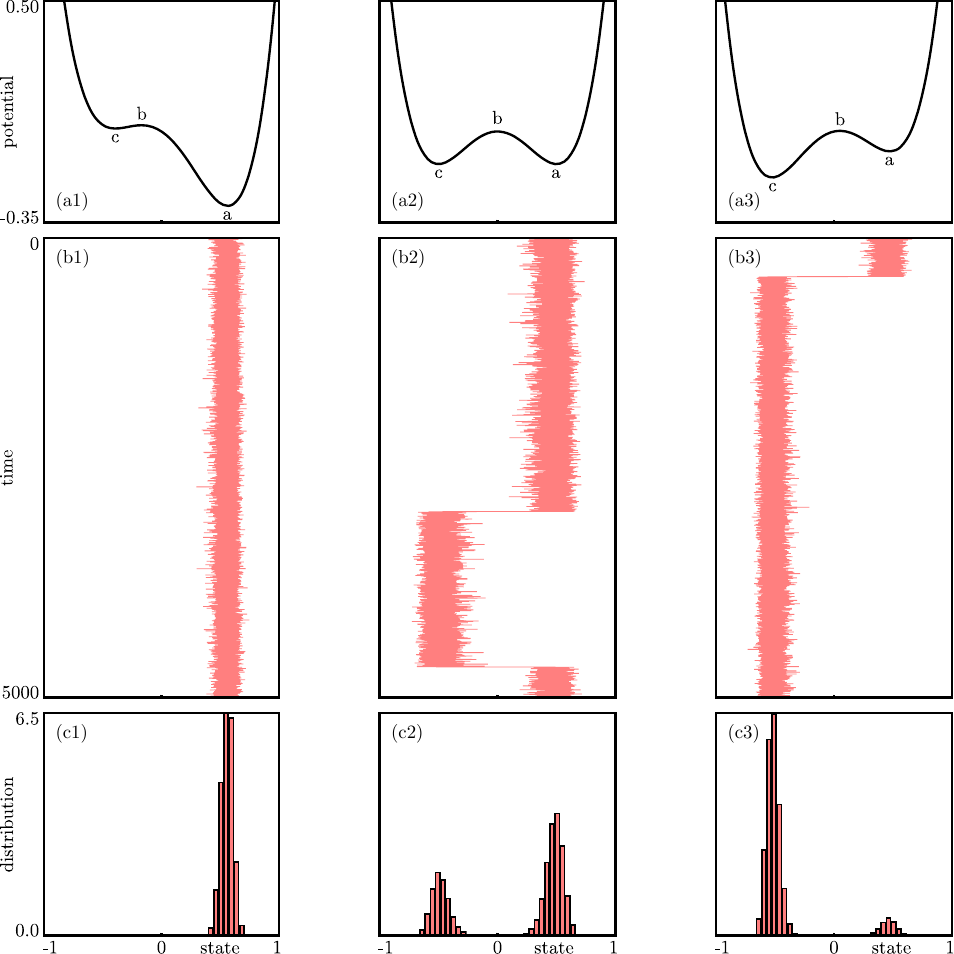}
        \caption{Different configurations of a generic double-well potential leading to qualitatively different trajectories of Eq.~\eqref{eq:langevin}.
                 Panels (a1-a3) show the shape of the confining potential with wells $a$ and $c$ separated by a hilltop $b$; panels (b1-b3) display a representative sample path of particle starting from the rightmost well $a$; panels (c1-c3) represent the sample distribution of the above timeseries.
                 In panels (a1-c1) the particle is confined by a high potential barrier and does not escape its basin of attraction within the simulation time.
                 In panels (a2-c2) the system is at a Maxwell point (both wells are equiprobable at stationarity) and the particle visits both wells during the simulation (the histogram is bimodal and becomes symmetric for a sufficiently long timeseries).
                 In panels (a3-c3) the starting well $a$ is shallower and flatter than the alternative one $c$: as such the particle escapes it early, never to return before the simulation ends.
                 The solutions of the Langevin equations of motions for each configuration have been computed using the DifferentialEquations.jl package \cite{Rackauckas2017}.}
        \label{fig:double_well_potential}
\end{figure}

For simplicity of exposition, let us consider the $1-$dimensional case of Eq.~\eqref{eq:langevin} for which we are always guaranteed to have a potential $V$ whose gradient reconstructs the dynamics $f$ in Eq.~\eqref{eq:fast_slow_sde}.
Examples of potentials $V$ are shown in \fref{fig:double_well_potential} (a1-a3).
By restricting ourselves to model single-parameter SN critical transitions (i.e. the parameter space is $\mathbb{R}$), we will consider the case of a bistable system with at least one stable equilibrium for any $\mu$.
We note that higher-dimensional systems will generically undergo a SN bifurcation, and that it is in principle possible to change to a coordinate frame (aligning with the center manifold) in which \eqref{eq:langevin} describes the relevant bifurcating dynamics. 
In practice the appropriate coordinate transformation is unknown, and an open problem is to determine optimal observables of high-dimensional systems \cite{Lohmann2025b}.

The bistability region of our system is characterized by a subset $\mu_1 < \mu < \mu_2$ of the parameter space for which two stable equilibria coexist alongside a single unstable equilibrium.
The bistability region is bounded by two SN bifurcations at $\mu=\mu_1$ and $\mu=\mu_2$.
The resulting potential features two local minima or \textit{wells} (labelled ``a'' and ``c'' in \fref{fig:double_well_potential}) separated by a barrier of one local maximum or \textit{hilltop} (labelled ``b'' in \fref{fig:double_well_potential}).
Different configurations of the shape of the potential describe qualitatively different particle dynamics as represented in \fref{fig:double_well_potential} (b1-b3), which corresponds to the potentials $V$ in \fref{fig:double_well_potential} (a1-a3).
In statistical mechanics, Kramer's formula quantifies the rate at which particles escape across potential barriers when subject to Langevin's equation of motion \eqref{eq:langevin}.
Kramer's formula is sometimes expressed as the reciprocal of the escape time, usually referred to as the (mean) first exit time.
Different derivations of the formula are available in the literature; we thus redirect the readers to \cite{Berglund2013} for a concise overview.
It must be noted that, while limitations on Kramer's law vailidity exist (see \cite{Berglund2013}), our assumptions of stationary, $1-$dimensional dynamics subject to small, additive white noise fall within the applicability cases.
 
Given a double-well potential $V$ with desirable local minimum (stable equilibrium) $x=a$ and local maximum (unstable equilibrium) $x=b$, the rate at which particles, starting at $x=a$ and subject to white noise of intensity $\sqrt{2D}$, will escape over the saddle $x=b$ is given by
 
\begin{equation}\label{eq:kramer_escape_formula}
        R_{a\to b} = \frac{\sqrt{|V''(b)|V''(a)}}{2\pi}\exp\bigg(-\frac{\Delta V_{a\to b}}{D}\bigg)\,,
\end{equation}
where $\Delta V_{a\to b} = V(b) - V(a)$ quantifies the potential barrier that the Langevin overdamped particle has to overcome in order to escape, and $V''(a),V''(b)$ quantify the local curvature of $V$ around the well and the hilltop respectively.
\subsection{An interpretable EWS}\label{subsec:modified_escape_ews}
Among the many EWS that have been proposed in the literature, the increase in variance has been the most popular.
Its success primarly stems from a relative simplicity in computation and a relative robustness in detecting CSD.
The major drawback of using the increase in variance, as an indication of B-tipping, is its lack of interpretability.
As discussed in the introduction, numerical values of the variance of a sample fail to provide useuful interpretation if one does not have access to a baseline value for the system when at a safe distance from B-tipping.

For the rest of this section we show that a modified version of Kramer's escape formula solves this issue by providing an interpretable EWS.
Consider a generic SN bifurcation in $d-$dimensions; if the system is conservative and the noise level is small, then \eqref{eq:kramer_escape_formula} holds when the state of the system is far enough from the bifurcation.
In this regime, the two wells of the potential are well separated and thus the local minimum (stable equilibrium) $x=a$ is much deeper than the local maximum (unstable equilibrium) $x=b$.
In other words, far from the bifurcation the potential barrier $\Delta V = V(b) - V(a) \gg 0$ (we dropped the transition notation $a\to b$).
Conversely, as the SN bifurcation is approached, the stable equilibrium $x=a$ and the unstable equilibrium $x=b$ of the system become closer, resulting in the potential barrier reducing, which implies that $\Delta V\to0$.

Based on the above, we introduce our modified escape EWS
\begin{equation}\label{eq:escape_ews}
        \exp(-\Delta V) \approx 
   \begin{cases}
       0\,, \quad |\mu - \mu_{c}| \gg 0\,,\\
       1\,, \quad \mu \to \mu_{c}\,. 
   \end{cases}
\end{equation}
where $\mu=\mu_c$ denotes the parameter value at which a SN bifurcation between the stable state $x=a$ and the unstable saddle $x=b$ occurs.
Notice that \eqref{eq:escape_ews} is obtained from \eqref{eq:kramer_escape_formula} by removing the curvature-dependent prefactor $\sqrt{|V''(b)|V''(a)}/2\pi$ and setting $D=1$ in the argument of the exponential.
Removing the prefactor from Kramer's escape rate \eqref{eq:kramer_escape_formula} guarantees that the modified formulation \eqref{eq:escape_ews} increases monotonically as B-tipping is approached (whereas with the prefactor in \eqref{eq:kramer_escape_formula} included, as the bifurcation is approached, Kramer's formula will initially increase away from zero, before reaching a maximum and then decreasing back to zero at the bifurcation).
Setting $D=1$ guarantees that the modified formulation is not late in detecting the proximity of B-tipping when the underlying noisy fluctuations are small, i.e. when $D\to0$  (see \aref{appendix_a} for a detailed discussion).
Furthermore, the true value of $D$ is generally not known for real systems and, while estimation is possible it inflates the dimensionality of the problem without providing significant advantages to the underlying EWS, as shown in \fref{fig:kramers_formula}(a) in \aref{appendix_a}.

The boundedness of \eqref{eq:escape_ews} in $(0,1]$ solves the interpretability problem we outlined above.
As opposed to the variance (and higher-order moments), a numerical value of the modified escape formula provides direct interpretation of the likelihood of a B-tipping event without the need for a baseline value to refer to.
If a system is far from a dangerous transition (i.e. $|\mu-\mu_c|\gg0\implies\Delta V\gg0$), then the modified escape EWS will be close to $0$; conversely, a state that is approaching said transition (i.e. $\mu\to\mu_c\implies\Delta V\to0$) results in values of the same quantity to be close to $1$.

\begin{figure}
	      \centering
	      \includegraphics[keepaspectratio, width=\textwidth]{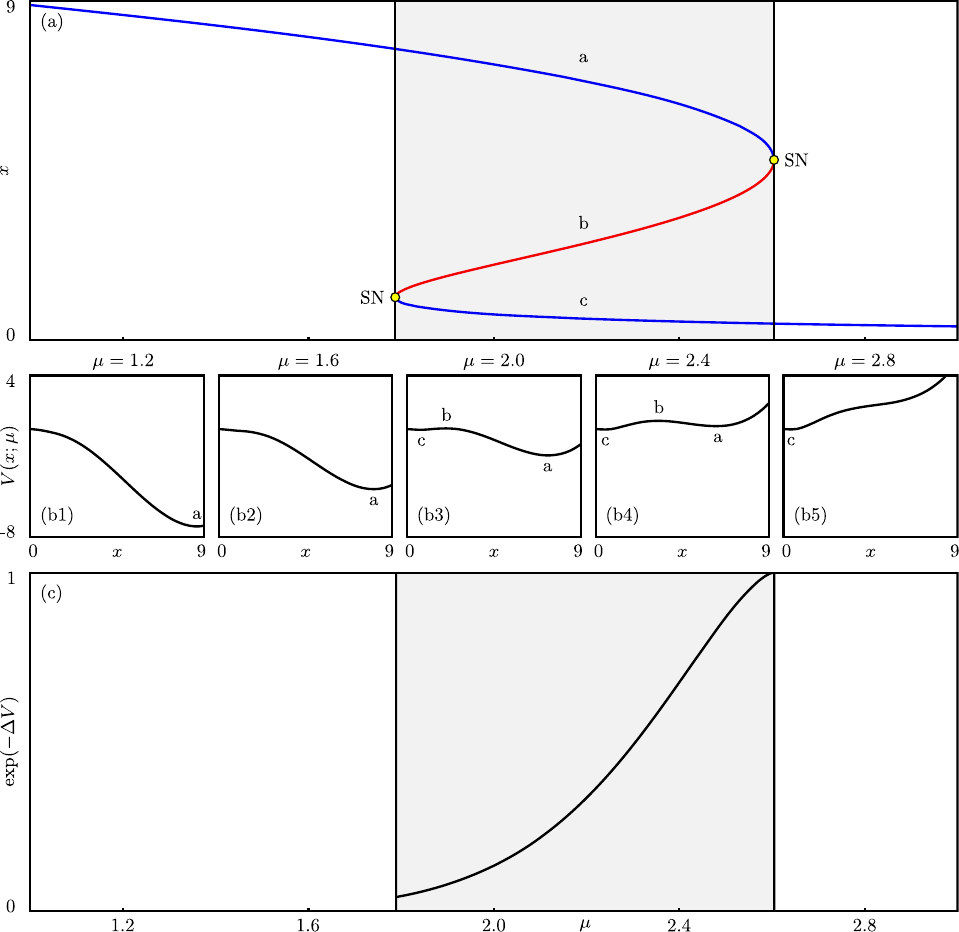}
        \caption{(Panel (a)) Bifurcation diagram of a model of ecosystem collapse with the state $x$ (vegetation biomass) degrading as the forcing parameter $\mu$ (harvesting rate) is increased \cite{May1977}; branches of stable equilibria ($x=a$ and $x=c$) are in blue while the unstable branch ($x=b$) is in red; the region of bistability, bounded by two SN bifurcations (yellow circles), is shaded in gray.
        (Panels (b1-b5)) Different configurations of the potential $V(x;\mu) = -\frac{x^2}{2} + \frac{x^3}{30} + \mu\big(x - \arctan(x)\big)$ at different parameter values (ticks of the horizontal axes of panels (a) and (c)). 
        (Panel (c)) Modified escape EWS \eqref{eq:escape_ews} computed out of the potential $V(x;\mu)$ between the stable state $x=a$ and the saddle $x=b$ within the bistability region.
        The bifurcation diagram has been computed using the BifurcationKit.jl package \cite{Veltz2020}.}
        \label{fig:escape_ews}
\end{figure}

In \fref{fig:escape_ews} we depict the monotonic increase of \eqref{eq:escape_ews} from $0$ to $1$ for a bistable model of ecosystem collapse under increasing harvesting \cite{May1977}.
The region of bistability (depicted by the gray area in \fref{fig:escape_ews} (a,c)) is bounded by two SN bifurcations at $\mu_1$ and $\mu_2>\mu_1$.
The earlier SN bifurcation at $\mu=\mu_1$ involves the stable equilibrium $x=c$ and the unstable saddle $x=b$ while the later SN bifurcation at $\mu=\mu_2$ involves the stable equilibrium $x=a$ and the same saddle $x=b$.
This bifurcation sequence is reflected in the shape of the potential $V$ shown at a number of parameter slices in \fref{fig:escape_ews} (b1-b5), wherein the local maximum $x=b$ and the local minima $x=a$ and $x=c$ correspond to the unstable and stable equilibria in \fref{fig:escape_ews} (a), respectively.
Outside the shaded region (i.e. $\mu<\mu_1$ or $\mu>\mu_2$, e.g. \fref{fig:escape_ews} (b1,b2,b5)) only one stable equilibrium exists (either $x=a$ or $x=c$), the potential $V$ has one global minimum and no local maxima and thus the modified escape EWS is undefined.
Conversely, inside the shaded region (i.e. $\mu_1<\mu<\mu_2$, e.g. \fref{fig:escape_ews} (b3,b4)) the potential $V$ has two minima ($x=a$ and $x=c$) and one local maximum ($x=b$), meaning that the escape EWS can be computed using \eqref{eq:escape_ews} (e.g. \fref{fig:escape_ews} (c)).

Notice that the modified escape EWS \eqref{eq:escape_ews} increases monotonically with $\mu$, with its minimum ($\exp(-\Delta V)\approx0$) being at $\mu=\mu_1$ and its maximum ($\exp(-\Delta V)=1$) being at $\mu=\mu_2$.
Therefore, by identifing the subset of the parameter space $\mathbb{R}$ such that $\Delta V = V(b) - V(a)$ has range $[0,+\infty)$, we associate to the modified escape EWS the intepretation of a probability measure, quantifying the likelihood of a B-tipping event to occur.
\section{Maximum likelihood estimation of the potential}\label{sec:mle}
\fref{fig:escape_ews} shows that \eqref{eq:escape_ews} can be used as an interpretable EWS whenever the potential $V$ of a gradient system \eqref{eq:langevin} is known which, in most cases, it is not.
As mentioned above, the popularity of canonical EWS, such as the increase in timeseries variance, is in part due to their simplicity in calculating them from available data.
The increase in variance for example can be readily detected from the timeseries of a single realization of Eq.~\eqref{eq:langevin}.
In order for our proposed, interpretable EWS to be applicable in the real-world we must be able to compute it from timeseries data without relying on prior knowledge of the system that generates it.
We will now show that, under certain constraints, a similar procedure of extracting EWS from timeseries data can be used for the modified escape rate \eqref{eq:escape_ews}.

Since the modified escape rate depends on the potential barrier $\Delta V$ between the observed steady-state $x=a$ and an unknown unstable equilibrium $x=b$, our aim for this section is to estimate such quantity from a single timeseries. 
To do so we must be able to approximate the local shape of the potential $V(x;\mu)$ between $x=a$ and $x=b$.
As we are trying to infer the likelihood of a SN B-tipping, the minimum model that we can fit to a timeseries is that of a cubic polynomial $V(x;\theta) = \sum_{k=1}^{3}\theta_{k}x^{k}$ with real coefficients $\theta=(\theta_{1},\theta_{2},\theta_{3})$.
While different choices are available for our inference problem, a cubic polynomial is minimal in the sense that it allows for the presence of one local minimum $x=a$ and one local maximum $x=b$ (from which $\Delta V = V(b) - V(a)$ can be computed) without inflating the dimensionality of the search space.

The inference of the polynomial coefficients of the potential $V$ is, at its core, a parameter estimation problem of the deterministic drift of a SDE given a sequence of its realizations.
This is a well-established problem in statistical inference \cite[Ch. 11, p. 234]{Sarkka2019}; in the following we briefly formulate the setup for such problem and redirect interested readers to \cite{Nielsen2000, Sorensen2004, Craigmile2022} for reviews of these techniques.
\subsection{Problem formulation}\label{subsec:formulation}
A standard technique to infer the parameters of a model SDE, using its solution, is maximum likelihood estimation (MLE).
Suppose our sample $\mathcal{X}=\{X_{n}\}_{n=0,1,\dots,N}$ is made of $N+1$ observations of the stationary state solving Eq.~\eqref{eq:langevin} from time $t=0$ to time $t=T$.
Suppose also that each observation is sampled at homogeneous frequency of time intervals $\Delta t = T/N$.
This It\^{o} process is a Markov process \cite[Ch.7, p. 111]{Oksendal2000} which means that the transition probability density to a future state $X_{t+1}$, at any time $t$, solely depends on the current state $X_{t}$ rather than the entire history $\mathcal{X}$.
This results in the joint distribution of the sample to be written in a factorized form of the (marginal) transition probability densities
\begin{equation}\label{eq:joint_distribution}
        p(\mathcal{X}) = p(X_{0})\prod_{n=1}^{N-1}p(X_{n+1}|X_{n})\,.
\end{equation}
From \eqref{eq:joint_distribution} we can define the likelihood of the sample to be generated by a parametrized transition density 
\begin{equation}\label{eq:likelihood_function}
        \mathcal{L}(\theta|\mathcal{X}) := p(X_{0}|\theta)\prod_{n=1}^{N-1}p(X_{n+1}|X_{n},\theta)\,,
\end{equation}
with parameters $\theta\in\mathbb{R}^p$.
By maximizing \eqref{eq:likelihood_function} one finds the most likely set of parameters $\theta^{*}$, for a family of transition densities $p(x|\theta)$, to have generated the observations $\mathcal{X} = \{X_{n}\}_{n=0,\dots,N}$.
Notably one reformulates this problem in the minimization of the inverse log-likelihood instead
\begin{equation}\label{eq:inv_log_likelihood}
        \ell(\theta|\mathcal{X}) = -\log\mathcal{L}(\theta|\mathcal{X}) = -\log p(X_{0}|\theta) - \sum_{n=1}^{N-1}\log p(X_{n+1}|X_{n},\theta) \approx -\sum_{n=1}^{N-1}\log p(X_{n+1}|X_{n},\theta) \,,
\end{equation}
as this does not change the stationary point in the search space
\begin{equation}\label{eq:mle}
        \theta^{*}(\mathcal{X}) = \argmax_{\theta\in \mathbb{R}^{p}}\mathcal{L}(\theta|\mathcal{X}) = \argmin_{\theta\in \mathbb{R}^{p}}\ell(\theta|\mathcal{X})\,.
\end{equation}
Notice that in \eqref{eq:inv_log_likelihood} we approximate the inverse log-likelihhod by ignoring the conditional probability of the deterministic initial condition $X_{0}$ as this term is negligible compared to the sum of the log-likelihoods of the stochastic observations. 
The optimization problem \eqref{eq:inv_log_likelihood} requires the knowledge of the transition density $p(X_{n+1}|X_{n},\theta)$, which, for a SDE of the form \eqref{eq:langevin}, is goverened by the Fokker-Plank equation (FPE)
\begin{equation}\label{eq:fpe}
        \partial_{t}p(x,t|\theta) = \partial_{x}(V'(x;\theta)\,p(x,t|\theta)) + D\,\partial_{xx}^{2}p(x,t|\theta)\,.
\end{equation}
Solutions of the FPE are rarely available in closed form and thus the optimization \eqref{eq:inv_log_likelihood} often relies on numerical or analytical approximations. 
\subsection{Time-invariant approximation: the Ornstein-Uhlenbeck process}\label{subsec:oup}
In the general time-varying regime, solutions of the FPE \eqref{eq:fpe} are usually found numerically.
However, in the specific case of slow-fast systems of the form \eqref{eq:fast_slow_sde}, the computation of the transition density can be simplified. 
Specifically, in the singular limit ($\varepsilon\to0$) the (linear) ramp of the forcing parameter is much slower than the timescale of the evolution for the state variable. 
This means that the time-varying solution of Eq.~\eqref{eq:fpe} can be reasonably approximated with a time-invariant density $p(x|\theta)\approx p(x,t|\theta)$ by neglecting the (partial) time derivative $\partial_{t}p\approx0$.
This ansatz reduces the parabolic partial differential equation in time and space \eqref{eq:fpe} to a second-order, homogeneous ordinary differential equation in $p(x|\theta)$ with non-constant coefficients
\begin{equation}\label{eq:stationary_fpe}
        0 = \partial_{x}\big(V'(x;\theta)\,p(x|\theta)\big) + D\,\partial_{xx}^{2}p(x|\theta) = V''(x;\theta)\,p(x|\theta) + V'(x;\theta)\,p'(x|\theta) + D\,p''(x|\theta)\,.
\end{equation}
The solution of the above (under suitable choices of boundary conditions at $x\to\pm\infty$) yields a Boltzmann-like probability density function \cite[Ch. 5, p. 98]{Risken1989}
\begin{equation}\label{eq:invariant_pdf}
        p(x|\theta) = \frac{1}{Z(\theta)}\exp\bigg(-\frac{V(x;\theta)}{D}\bigg)\,,
\end{equation}
with the normalization constant $Z(\theta)$ playing the role of the partition function of the energy of the system.
Thus, in principle, one can use the time-invariant solution \eqref{eq:invariant_pdf} of the FPE as the transition density in the optimization problem \eqref{eq:mle}.
 
\begin{figure}[t]
    \centering 
    \includegraphics[keepaspectratio, width=\textwidth]{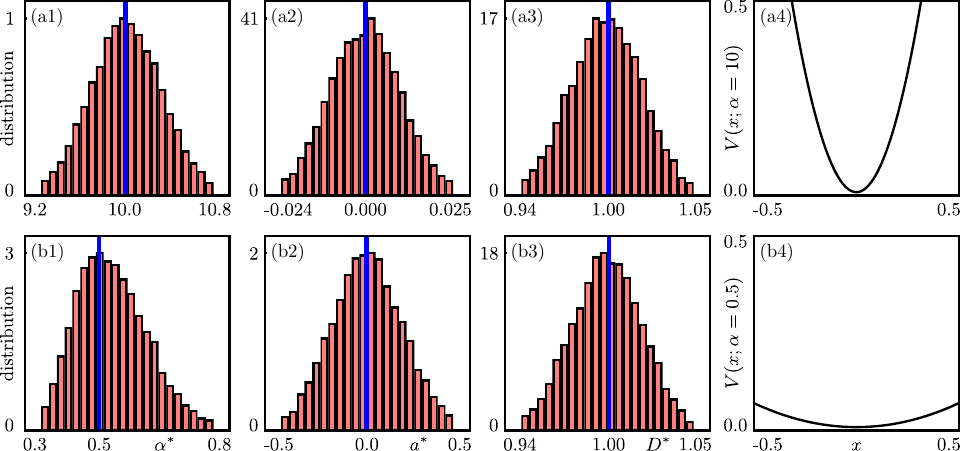}
    \caption{Empirical distribution of the MLE of the parameters $\theta^* = (\alpha^*, a^*, D^*)$ of an ensemble of $10,000$ sample paths solving \eqref{eq:oup} for two different configurations of the potential $V(x) = \alpha(a - x)^{2}$: a deep, strongly confining well with $\alpha = 10$ (panel (a4)) and; a shallow well with $\alpha = 0.5$ (panel (b4)).
            In each of the panels (a1-a3) and (b1-b3) the true value of the parameter is reported (blue line).
            Note, as $\alpha\to0$ the distribution of $\alpha^*$ becomes skewed despite the symmetry in the confining potential (panel (b1)).
            The solutions of the MLE problem have been computed using the pymle package \cite{Kirby2025}.}
    \label{fig:oup_mle}
\end{figure}
 
As an example of applicability of this workflow, let us consider the well-known case of an Ornstein-Uhlenbeck process (OUP) of the form
\begin{equation}\label{eq:oup}
        dX_{t} = \alpha(a - X_{t})dt + \sqrt{2D}dW_{t},
\end{equation}
where $\alpha,a> 0$.

The illustrative value of considering the OUP is given by the fact that it approximates the solution of any nonlinear (stationary) SDE around an equilibrium $x=a$ of the drift term $f$ in \eqref{eq:fast_slow_sde}.
More specifically, if \eqref{eq:fast_slow_sde} is not ramped (i.e. $\varepsilon=0$), the sample path around a stable equilibrium $x=a$ of the nonlinear vector field $f(x) = -V'(x)$ is well approximated by the solution of \eqref{eq:oup} with $\alpha = f'(x=a) = -V''(x=a) > 0$.
By denoting $\theta=(\alpha,a,D)$ to be our parameter vector, we setup a MLE problem in which the (limiting) Gaussian distribtion of the OUP \cite[Ch. 3, p. 74]{Gardiner2009}
\begin{equation}\label{eq:oup_density}
        p(x|\theta) = \sqrt{\frac{\alpha}{\pi D}}\exp\bigg(-\frac{\alpha(a - x)^{2}}{D}\bigg)\,,
\end{equation}
is of Boltzmann-like form \eqref{eq:invariant_pdf} with partition function $Z(\theta) = \sqrt{\pi D/\alpha}$ and quadratic potential $V(x) = \alpha(a - x)^{2}$.
For large enough $\alpha$, the solution of this MLE problem has been shown to be asymptotically normal \cite{Sahalia2002}, as demonstrated in \fref{fig:oup_mle} (a1-a3).

When $\alpha\to0$ however, the distribution of $\alpha^*$ becomes skewed towards larger values (see \fref{fig:oup_mle} (b1)). 
On top of this, the standard deviation of the distribution of $a^*$ increases (see \fref{fig:oup_mle} (b2)), yielding an overall less accurate estimation.

\begin{figure}[!tp]
    \centering 
    \includegraphics[keepaspectratio, width=\textwidth]{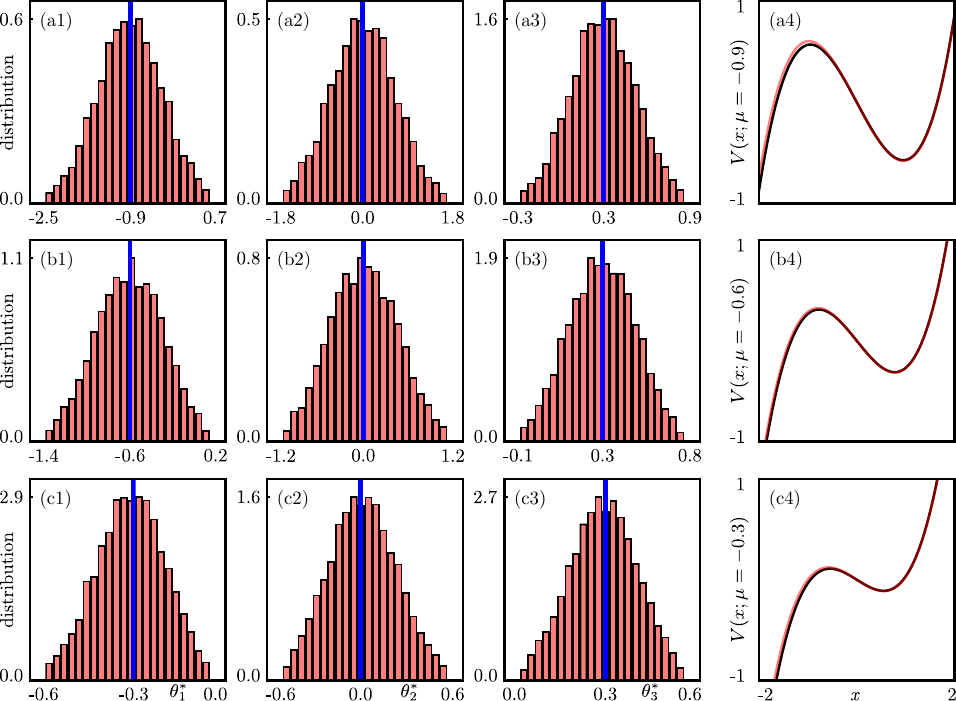}
    \caption{Empirical distribution of the solutions $\theta^{*} = (\theta_{1}^{*},\,\theta_{2}^{*},\,\theta_{3}^{*})$ of the LLS problem \eqref{eq:lls_solution} for an ensemble of $5,000$ sample paths of a SN normal form $f(x;\mu) = -\mu - x^2$.
            The problem is solved for different configurations of the potential $V(x;\mu) = \mu x + \frac{x^3}{3}$; at $\mu = -0.9$ (panels (a1-a4)), $\mu=-0.6$ (b1-b4) and $\mu=-0.3$ (c1-c4).
            The value of the polynomial coefficients of the true potential are reported in each of the panels (a1-a3), (b1-b3) and (c1-c3) (blue line).
            The ensemble mean values of the LLS solutions $\hat{\theta}^{*} = (\hat{\theta}_{1}^{*},\,\hat{\theta}_{2}^{*},\,\hat{\theta}_{3}^{*})$ are used to plot an ensemble mean reconstruction (red curve in (a4-c4)) of the potential $\hat{V}^{*}(x) := V(x;\hat{\theta}^{*}) = \sum_{k=1}^{3}\hat{\theta}_{k}^{*}x^{k}$ and compared with the ground truth (black curve in (a4-c4)).}
    \label{fig:lls_solutions}
\end{figure}

Since we can interpret $\alpha=f'(x=a)$ to be the leading eigenvalue of the linearized vector field $f$ around the stable equilibrium $x=a$, we conclude that the MLE of a nonlinear SDE using the OUP transition density decreases in accuracy as the B-tipping event is approached.

\subsection{Euler-Maruyama least-squares regression}\label{subsec:euler}
A desirable property of an EWS is that its prediction becomes more accurate and reliable when the monitored system is in a proximity of the critical transition.
Conversly in fact, an EWS that becomes less accurate when the system is close to B-tipping will systematically fail to detect the incoming critical transition by providing false negative outcomes, which has, arguably, a more catastrophic implication than providing false positive outcomes.
Estimating the likelihood of B-tipping using the transition density of the stationary OUP however achieves the opposite of this, as its accuracy decreases the closer the sample path is to the SN bifurcation.
In this section we remedy and overcome this limitation by using an alternative approach to MLE.
To model B-tipping out of a sample path of \eqref{eq:fast_slow_sde} we fit a cubic polynomial $V(x;\theta) = \sum_{k=1}^{3}\theta_{k}x^{k}$ to its potential.
Notice that, in order for the stationary distribution $p(x)$ to be a solution of the FPE \eqref{eq:fpe} with cubic $V$ (quadratic $V'$), we need to truncate the domain of the probability density function to be either $[b,+\infty)$ or $(-\infty,b]$, with $x=b$ being the local maximum of $V$, depending on the sign of the cubic coefficient $\theta_{3}$.
This is equivalent to imposing reflecting boundary conditions at $x=b$ to the FPE.

A number of computational methods to approximate the stationary density $p(x)$ are available \cite{Pedersen1995, Jensen2002, Sahalia2002}.
Conversely, other techniques have been devised that do not rely on the numerical approximation of the transition density such as the generalized method of moments (GMM) \cite{Hansen1982}, Markov chain Monte-Carlo (MCMC) algorithms \cite{Roberts2001,Alexander2005} and quasi-maximum likelihood regressions (QML) \cite{Castelle1986, Yoshida1992}.
Here we focus on the latter approach and in particular we describe one QML technique that reduces the MLE problem into a linear-least squares (LLS) regression.

As usual, we assume our sample path $\mathcal{X} = \{X_{0},\dots,X_{N}\}$ to be a solution of \eqref{eq:langevin} under the stationarity restriction.
We can then approximate the entire sample path using the Euler-Maruyama time-stepping scheme.
This defines a stochatic iterated map of the form
\begin{equation}\label{eq:euler_maruyama}
        X_{n+1} = X_{n} -V'(X_{n};\theta)\Delta t + \sqrt{2D\Delta t}\,\xi\,,\quad \xi\sim\mathcal{N}(0,1)\,,
\end{equation}
where $\theta=(\theta_1,\theta_2,\theta_3)\in\mathbb{R}^{3}$ is the vector of polynomial coefficients parameterizing the cubic potential.
We rearrange Eq.~\eqref{eq:euler_maruyama} to define the first-order time increments
\begin{equation}\label{eq:increments}
        Y_{n} := \frac{X_{n+1} - X_{n}}{\Delta t} = -V'(X_{n};\theta) + \sqrt{\frac{2D}{\Delta t}}\,\xi = -V'(X_{n};\theta) + \Xi\,,\quad\Xi\sim\mathcal{N}\bigg(0,\frac{2D}{\Delta t}\bigg)\,.
\end{equation}
Notice now how \eqref{eq:increments} defines a regression to a model $V'$ of some data $\{X_{n}\}_{n=0,\dots,N-1}$ given a set of $N$ observations $\{Y_{n}\}_{n=0,\dots,N-1}$ perturbed by some Gaussian noise $\Xi$.
This sets up a parameter estimation problem in which one wishes to find an optimal parametrization $\theta^*$ of the potential such that the residual between the model predictions $V'(X_{n};\theta^*)$ and the observations $Y_{n}$ is minimized in some sense.
As anticipated above, we formulate a LLS problem by interpreting the noise $\Xi$ to be an unresolved fluctuation error in the observed data.
We thus minimize the residual vector
\begin{equation}\label{eq:residual}
        r(\theta) = \big(Y_{n} + V'(X_{n};\theta)\big)_{n=1,\dots,N-1}\,,
\end{equation}
between the observations $Y_n$ and the fitting model $V'(X_n,\theta)$ in the least-squares sense.
That is, we solve the minimization problem
\begin{equation}\label{eq:lls_problem}
        \theta^{*} = \argmin_{\theta\in \mathbb{R}^{3}}||r(\theta)||_{2}\,,
\end{equation}
\begin{algorithm}[!t]
        \DontPrintSemicolon
        \SetAlgoNoLine
        \SetAlgoNlRelativeSize{0}
        \SetNlSty{text}{\Indp}{.}
        \SetKwInput{Input}{Input}
        \SetKwInput{Output}{Output}
        \BlankLine 
        \Input{sample $\mathcal{X}=\{X_0,X_1,\dots,X_N\}$ (a single timeseries)}
        \Output{estimated proximity to B-tipping}
        \BlankLine 
        \hrule
        \BlankLine
        \nl Compute quasi-stationary residuals $\overline{\mathcal{X}}$ using a detrending technique.
        \BlankLine
        \nl Assemble the LLS problem \eqref{eq:increments}\;
        \Indp\For{$(X_n,X_{n+1})\in \overline{\mathcal{X}}$}{
                \vspace{-10pt}
                \begin{align*}
                        y_n &= X_{n+1} - X_n\,,\\
                        A_n &= \big(X_n,\;X_n^2,\;X_n^3\big)\,.
                \end{align*}
                \vspace{-20pt}
        }
        \Indm\BlankLine
        \nl Solve the normal equation \eqref{eq:lls_solution}
        \begin{equation*}
                \theta^* = (\theta^*_1, \theta^*_2, \theta^*_3) = (A^TA)^{-1}A^Ty\,.
        \end{equation*}
        \BlankLine
        \nl Compute the local minimum and local maximum of $V^*(x)=\sum_{k=1}^{3}\theta_k^*x^k$
        \begin{align*}
                x_{\text{min}} &= + \frac{1}{3\theta^*_3}\bigg(\sqrt{(\theta_2^*)^2 - 3\theta_1^*\theta_3^*} - \theta^*_2\bigg)\,,\\
                x_{\text{max}} &= - \frac{1}{3\theta^*_3}\bigg(\sqrt{(\theta_2^*)^2 - 3\theta_1^*\theta_3^*} + \theta^*_2\bigg)\,.
        \end{align*}
        \BlankLine
        \nl Compute the potential barrier 
        \begin{equation*}
                \Delta V^*=V^*(x_{\text{max}}) - V^*(x_{\text{min}})\,.
        \end{equation*}
        \BlankLine
        \nl Extract the modified escape EWS using Eq.~\eqref{eq:escape_ews}
        \begin{equation*}
                \text{EWS} = \exp(-\Delta V^*)\,.
        \end{equation*}
        \vspace{-10pt}
        \caption{Extraction of the modified escape EWS from timeseries data.}
        \label{alg:modified_ews}
\end{algorithm}
where $||\cdot||_2$ denotes the $L_2-$norm.
Notice that by fitting the minimal working model of B-tipping (i.e. we assume $V(x;\theta)$ is a cubic polynomial of coefficients $\theta = (\theta_{1},\, \theta_{2},\, \theta_{3})$), the residual \eqref{eq:residual} is linear in $\theta$.
This allows us to rewrite \eqref{eq:lls_problem} as a quadratic minimization problem
\begin{equation}\label{eq:quadratic_minimization_problem}
        \theta^{*} = \argmin_{\theta\in \mathbb{R}^{3}}||A\theta - y||_2\,,
\end{equation}
where $y\in \mathbb{R}^{N}$ denotes the observation data vector (i.e. $y_j = Y_j$), $A\in \mathbb{R}^{N\times3}$ is the model matrix evaluated on the sample $\mathcal{X}$ (i.e. $A_{jk} = X_{j}^{k}$) and $A\theta - y = r(\theta) \neq 0$.
The unique solution of \eqref{eq:quadratic_minimization_problem} satisfies the so-called normal equation
\begin{equation}\label{eq:lls_solution}
        \theta^{*} = (A^{T}A)^{-1}A^{T}y = A^{+}y\,,
\end{equation}
where $A^TA$ is a Gram matrix and $A^{+}\in\mathbb{R}^{3\times N}$ is the Moore-Penrose inverse of $A$.

With the estimation $\theta^*=(\theta^*_1, \theta^*_2, \theta^*_3)$ one can readily compute the modified escape EWS \eqref{eq:escape_ews} of the approximation $V^*(x):=V(x;\theta^*)$ using the sample path. We report the main steps of this procedure in Algorithm~\ref{alg:modified_ews}.
The solution of the normal equations \eqref{eq:lls_solution} is also asymptotically normal as shown in \fref{fig:lls_solutions} for the case of a SN normal form.

Notice that, contrary to the MLE using the stationary density of the OUP described in \sref{subsec:oup}, the numerical accuracy of \eqref{eq:lls_solution} improves as the SN bifurcation is approached.
This is shown in \fref{fig:lls_solutions} by a reduction of the standard deviation of the distribution (histograms) of the LLS solutions $\theta^*$ as $\mu\to0$.
We remark that this property is generic (i.e. it is observed for gradient systems other than the SN normal form) and, we argue, it could be intuitevely explained by critical slowing down (CSD).
As B-tipping is approached (i.e. $\mu\to0$) the leading eigenvalue of the linearized system approaches $0$ from below, which in turn entails that the stable equilibrium recovers perturbations of the steady-state with less \textit{``strength''}.
With the fluctuations around the stable equilibrium $x=a$ increasing in range, the state of the system will explore larger portions of the basin of attraction of such equilibrium.
As a consequence, the realizations $\mathcal{X} = \{X_0,\dots,X_N\}$ of the gradient system will capture the local shape of the confining potential in an increasingly better fashion as the system is driven closer to B-tipping, which translates in more accurate solutions of the LLS regression \eqref{eq:lls_solution}.

We conclude that, on top of its interpretability, the estimation of the escape EWS \eqref{eq:escape_ews}, as detailed in Algorithm \ref{alg:modified_ews}, becomes more trusthworty as B-tipping is approached, yielding a decreasing rate of false negatives.
For a more detailed and formal treatment of false negatives detection, and the temporal interpretation of proximity to B-tipping of CSD-based indicators, such as the one estimated by Algorithm~\ref{alg:modified_ews}, we refer interested readers to \cite{Ashwin2025}.

\section{Applications}\label{sec:applications}
We are now ready to extract the modified escape EWS \eqref{eq:escape_ews} from timeseries data by approximating the local shape of the confining potential using the likelihood estimator described in the \sref{subsec:euler}.
We will demonstrate the efficacy of the proposed modified escape EWS \eqref{eq:escape_ews} and its approximation method (Algorithm~\ref{alg:modified_ews}) using three test cases:
\begin{itemize}
     \item slowly ramped SN normal form with a single basin of attraction and computable (approximate) variance from the stationary solution of the FPE \eqref{eq:fpe};
     \item $1-$dimensional dynamics of vegetation biomass collapse in a regime of bistability under slowly increasing harvesting rate \cite{May1977};
     \item $1-$dimensional observable of a $2-$dimensional formulation of the thermohaline circulation in the North Atlantic driven by increasing freshwater forcing \cite{Stommel1961}.
\end{itemize}
In all of these tests, the modified escape EWS will be estimated by applying Algorithm~\ref{alg:modified_ews} to sliding windows of the slowly ramped solutions.
The windowed timeseries are then detrended using standard, parameter free techniques and the obtained quasi-stationary residuals retain similar statistical properties (see \sref{sec:analysis} for more details).
\subsection{Slowly ramped SN normal form}\label{subsec:SN_normal_form}
We follow the same setup outlined in \cite{Kuehn2011}.
A slow-fast framework of the form \eqref{eq:fast_slow_sde} is adopted for a SN normal norm
\begin{equation}
        \begin{aligned}\label{eq:slow_fast_SN}
                dX_{t} &= (-\mu-X_{t}^{2})dt + \sigma dW_{t}\,,\\
                d\mu &= \varepsilon dt\,, 
        \end{aligned}    
\end{equation}
where $\varepsilon\ll1$ is the timescale separation between the dynamics of the fast variable $X_{t}$ and that of the slow variable $\mu$.
We assume also that $0 < \sigma \ll \sqrt{\varepsilon} < 1$ so that solutions of \eqref{eq:slow_fast_SN} undergo B-tipping close to $\mu=\mu_{c}=0$ when ramped.
\begin{figure}[!tp]
    \centering 
    \includegraphics[keepaspectratio, width=\textwidth]{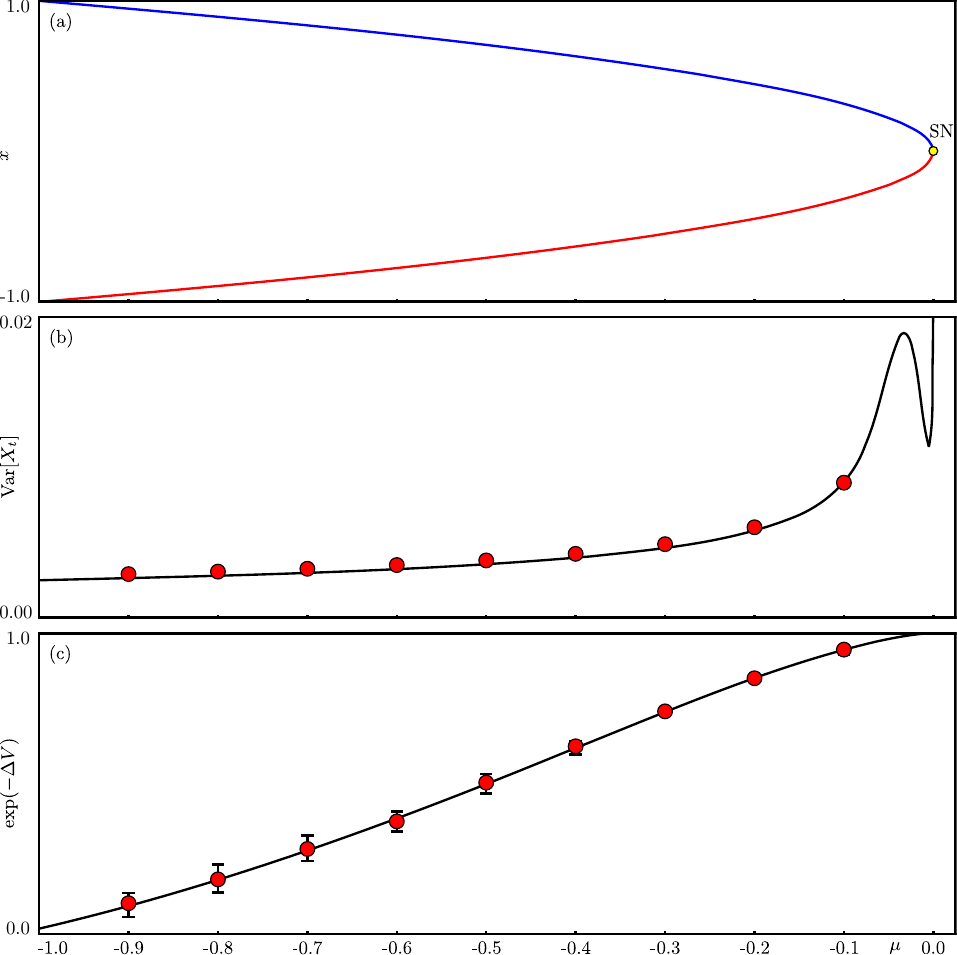}
    \caption{Comparison between the increase in variance and modified escape EWS for stationary solutions of a SN normal form \eqref{eq:slow_fast_SN} at different values of the bifurcation parameter $\mu\in(-1,0)$.
            (Panel (a)) Bifurcation diagram of a SN normal form with the blue curve representing the branch of the stable equilibrium $a = \sqrt{-\mu}$ while the red curve represents the unstable equilibrium branch $b = -\sqrt{-\mu}$; the SN bifurcation is depicted by the yellow circle at $\mu=0$.
            (Panels (b) and (c)) Comparison between the analytical (black curves) and empirical (red circles) variance and modified escape EWS.
            The analytical values have been computed using \eqref{eq:stationary_SN_variance} for the variance (numerical integration is performed using a Gauss-Kronrod quadrature rule \cite{Laurie1997}) and \eqref{eq:escape_ews} for the escape EWS; the empirical values have been computed by taking the median of an ensemble of $100$ sample paths at each parameter value.
    The error bars represent the interquartile range.}
    \label{fig:stationary_SN_ews}
\end{figure}
In the singular limit $\varepsilon\to0$ the solution of the time-invariant FPE gives a Boltzmann-like (stationary) probability density of the form
\begin{equation}\label{eq:stationary_SN_density}
        p(x;\mu) = \frac{1}{Z(\mu)}\exp\bigg(\frac{\mu x + \frac{x^{3}}{3}}{D}\bigg)\,.
\end{equation}
with $Z(\mu) = \int_{-\sqrt{-\mu}}^{+\infty}\exp\big(\frac{\mu x + \frac{x^{3}}{3}}{D}\big)dx$ being the normalization constant.
The evaluation of \eqref{eq:stationary_SN_density} has to be done numerically with $Z(\mu)$ being approximated by quadrature.
With this, one can also derive the numerical approximation of the analytical variance of the r.v. $X_{t}$ according to its definition
\begin{equation}\label{eq:stationary_SN_variance}
        \text{Var}[X_{t};\mu] = \int_{-\sqrt{-\mu}}^{+\infty}\big(x - \mathbb{E}[X_{t};\mu]\big)^{2}p(x;\mu)dx\,,
\end{equation}
which is again computed by quadrature.
In \eqref{eq:stationary_SN_variance}, $\mathbb{E}[X_{t};\mu] = \int_{-\sqrt{-\mu}}^{+\infty}x\,p(x;\mu)dx$ denotes the expected value of the r.v. $X_{t}$.

We first show that the modified escape EWS reconstructed from the MLE of the confining potential converges, in some sense, to its analytical value in the large number of samples limit.
\fref{fig:stationary_SN_ews} shows the result of a simulation of an ensemble of $100$ sample paths of Eq.~\eqref{eq:slow_fast_SN} for fixed values of the bifurcation parameter (i.e. with $\varepsilon=0$);
specifically, here and throughout the rest \sref{sec:applications}, we apply Algorithm~\ref{alg:modified_ews} to estimate the modified escape EWS for individual timeseries, and then analyse the statistics of an ensemble of timeseries realizations.
In doing so, we are able to say what is the typical and expected behavior of the EWS.
For each value of $\mu$ in the parameter sweep of the $(-1,0)$ interval and for each sample path in the ensemble, we solve the LLS problem \eqref{eq:lls_solution} of a fitting polynomial potential and compute the modified escape EWS \eqref{eq:escape_ews} (see Algorithm \ref{alg:modified_ews}).
This results in a distribution of modified escape EWS whose median we then compare to the ground truth, i.e. the analytical value of \eqref{eq:escape_ews} when evaluated using the true potential $V(x;\mu)=\mu x + \frac{x^3}{3}$. 
We do the same for the ensemble variance by computing \eqref{eq:stationary_SN_variance} via numerical approximation, using it as the \textit{``ground truth''}, and compare it against the median sample variance computed out of the solutions of \eqref{eq:slow_fast_SN}.
\begin{figure}[!tp]
    \centering 
    \includegraphics[keepaspectratio, width=\textwidth]{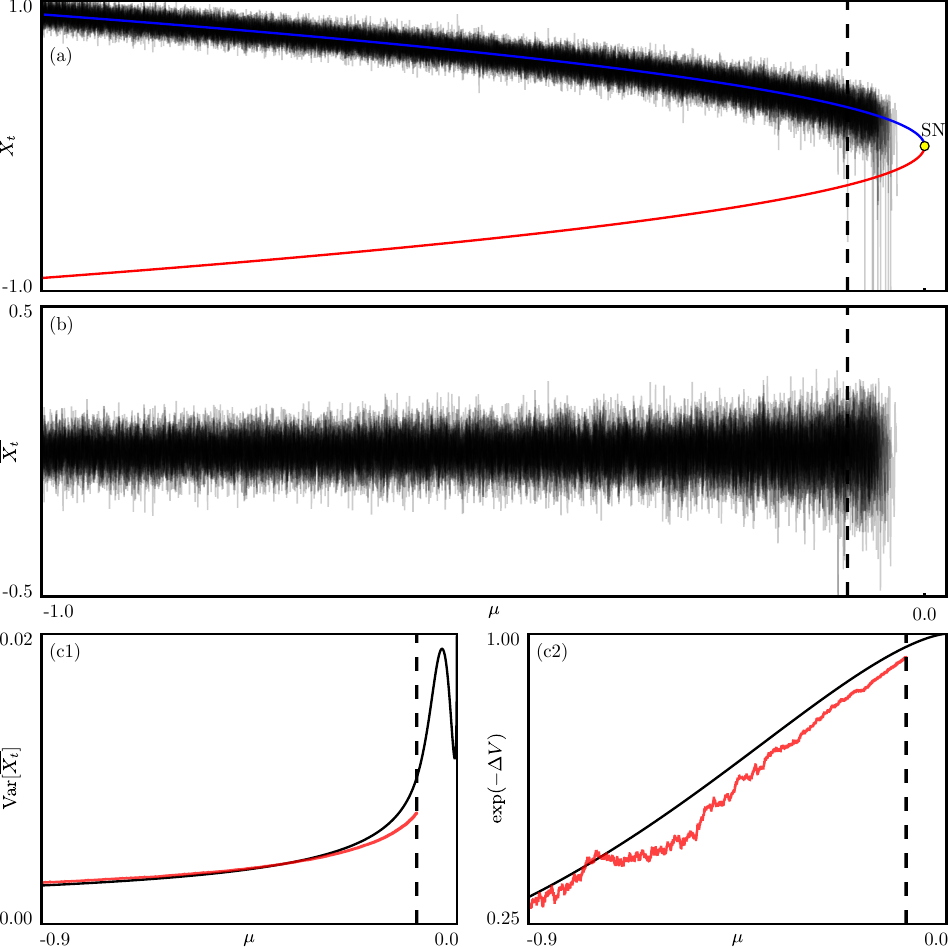}
    \caption{Comparison between the increase in variance and modified escape EWS for the solutions of the slow-fast formulation of a SN normal form \eqref{eq:slow_fast_SN}.
            (Panel (a)) An ensemble of $100$ sample paths (black), slowly forced with rate $\varepsilon=4\cdot10^{-5}$, follows the stable branch (blue) of the critical manifold $a(\mu)$; notice that all the sample paths tip across the unstable branch (red) before the SN bifurcation (yellow circle).
             We indicate the earliest tipping in the ensemble by the vertical dashed line. 
             (Panel (b)) Residuals $\overline{X_t}:=X_t - a(\mu)$ obtained by detrendeding the drift of the solutions $X_t$.
             The trend is computed using the pyEMD package \cite{Laszuk2017}. 
             (Panels (c1-c2)) Comparison between the analytic signal of the frozen system (black) against the ensemble averages of the same quantities computed (empirically) from the sample paths (red) on a sliding window of $10\%$ the length of the sample paths.}
    \label{fig:ramped_SN_EWS}
\end{figure}
 
The empirical results in \fref{fig:stationary_SN_ews} (c) show that, in the large number of samples limit, the modified EWS computed by MLE does converge to the true value in median.
This is not surprising given prior evidence of convergence in mean of the MLE coefficients to their true value as shown in \fref{fig:lls_solutions}.
We also emphasize that the accuracy of the MLE of the modified escape EWS improves as the system gets closer to the B-tipping as shown by the shrinking of the errorbars in \fref{fig:stationary_SN_ews} (c).
This is a highly desirable property of the proposed method as it minimizes the chances of false positives for dangerous critical transitions of the B-tipping type.
Furthermore, alongside interpretability, the modified escape EWS is also \textit{``earlier''} in the upward trend compared to the increase in variance which also shows non-monotonicity at the very onset of the B-tipping event.

Compared to the increase in variance however, the modified escape EWS estimated from the sample paths does have, on average, a larger standard deviation.
This becomes more evident when we allow the bifurcation parameter to vary slowly and apply the same methodology described above to detrended residuals of the solutions, as depicted in \fref{fig:ramped_SN_EWS}.
The source of this is embedded in the large dimensionality of the search space $\mathbb{R}^{3}\ni\theta$, correlations in the estimated coefficients and other well-known issues regarding the stability of the method (e.g. ill-conditioning of the Vandermonde matrix).
Regularization of the solutions using penalizers (e.g. Tikhonov regularizations) can alleviate numerical instabilities in the LLS method and while it reduces the spread of the solutions it also introduces biases that effectively bring the mean solution further away from the ground truth.
Furthermore, optimal choices of the penalizers are currently available a-posteriori (e.g. Morozov's discrepancy principle \cite{Nair2009}), meaning that prior knowledge of the underlying system is necessary in order for regularization to be effective.
See Appendix \ref{appendix_b} for a detailed discussion on the effects of regularization on LLS solutions.
\subsection{Ecosystem collapse}\label{subsec:ecosystem}
In the early years of EWS development, ecologists focused on a benchmark model of ecosystems collapse \cite{Scheffer2009,Dakos2012}.
First introduced by May in 1977 \cite{May1977}, this model of vegetation biomass degradation under increasing harvesting reads
\begin{equation}\label{eq:may_model}
        \dot{x} = rx\bigg(1-\frac{x}{k}\bigg) - \mu \frac{x^2}{x^2 + h^2}\,,
\end{equation}
where the (scalar) state variable $x$ represents the vegetation cover in an ecosystem, $r, K, h\in \mathbb{R}$ are fixed parameters modeling the biomass growth rate, carrying capacity and half-saturation respectively, and finally $\mu$ is the external harvesting forcing (bifurcation parameter).
\begin{figure}[t]
    \centering 
    \includegraphics[keepaspectratio, width=\textwidth]{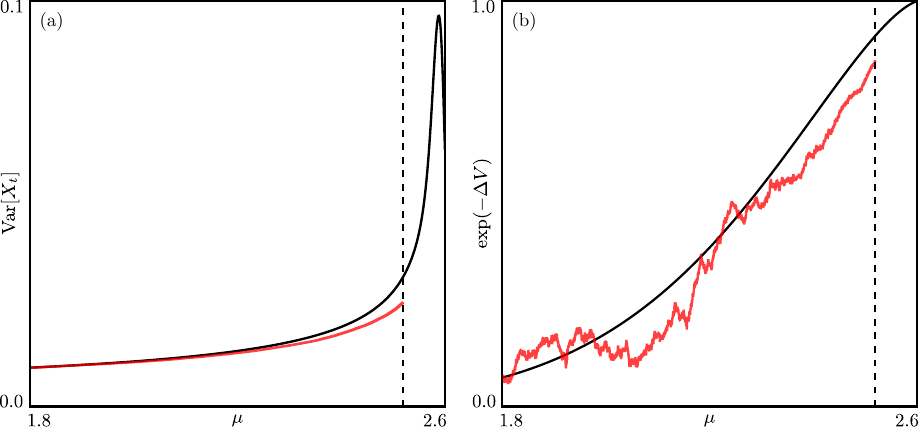}
    \caption{Comparison between the increase in variance and modified escape EWS applied to a slowly ramped model of ecosystem collapse \eqref{eq:may_model}, introduced in \cite{May1977}.
            We fix parameters $r=1$, $k=10$ and $h=1$ while we let $\mu$ to vary between $1.8$ and $2.6$ with rate $\varepsilon=5\cdot10^{-5}$; the bifurcation diagram of the stationary system is reported in \fref{fig:escape_ews}(a).
            A simulation of $100$ sample paths allows us to compute ensemble averaged variance and modified escape EWS (both in red) on a sliding window of $10\%$ the length of the solutions.
            In each panel the accuracy of empirical signals is compared to the approximations of the stationary system (both in black). 
            The earliest tipping in the solutions is shown by the dashed line.
            Notice how both variance and modified escape EWS show upward trends in the approach to the tipping point.
            The increase in variance however requires the knowledge of a baseline value which is different from that of other systems (e.g. see \fref{fig:ramped_SN_EWS}(c)) while the increase in escape likelihood shows a more reliable, direct interpretation as an EWS.}
    \label{fig:ramped_may_ews}
\end{figure}
Under a suitable choice of the fixed parameters, Eq.~\eqref{eq:may_model} shows a bistable region in $\mu$ where two stable states are accessible by the system.
A healthy stable equilibrium is characterized by high vegetation biomass; an alternative state of overexploitation coexists in the bistability region which corresponds to a low vegetation cover in the ecosystem.
The region of bistability is bounded by two SN bifurcations as depicted in \fref{fig:escape_ews} (a).
The potential of the system is non-polynomial
\begin{equation}\label{eq:may_potential}
        V(x;\mu) = -r\bigg(\frac{x^2}{2} - \frac{x^3}{3k}\bigg) + \mu\bigg(x - \sqrt{h}\arctan\bigg(\frac{x}{\sqrt{h}}\bigg)\bigg)\,,
\end{equation}
and the application of the MLE of the modified escape EWS to this system allows us to test its efficacy to more generic confining potentials.
We thus formulate a slow-fast problem of the form \eqref{eq:fast_slow_sde} with drift $f(x,\mu)$ equal to the right-hand side of \eqref{eq:may_model} while we let the harvesting forcing $\mu$ to slowly increase.
Similarly to the previous case, we can approximate the probability density of the process in the singular limit $\varepsilon\to0$ by looking at its stationary distribution, allowing us to compute the variance of the vegetation biomass numerically and compare it to an ensemble mean of $100$ sample paths solving the slow-fast formulation of \eqref{eq:may_model}.
Those sample paths will then form the dataset from which estimation of the cubic approximation of \eqref{eq:may_potential} is performed in order to extract the modified escape EWS.
The results of this experiment are depicted in \fref{fig:ramped_may_ews}.

What we notice in this case is that both signals show upward trends in the approach to the SN bifurcation.
The modified escape EWS provides a direct interpretation of its numerical values while the variance does not.
Furthermore we highlight that the increase in variance (\fref{fig:ramped_may_ews} (a)) shows a different range than the one measured for the SN normal form (\fref{fig:ramped_SN_EWS} (c1)).
The difference in ranges between different systems entails less confidence in interpreting the increase in variance as an early indicator of tipping, which is exacerbated in systems where data is scarce and the baseline value is unknown.
The issue is resolved by the (empirical) modified escape EWS which instead follows a similar range of values (\fref{fig:ramped_may_ews} (b)) as the one detected for the SN normal form (\fref{fig:ramped_SN_EWS} (c2)).
\subsection{Weakened circulation state in the North Atlantic}\label{subsec:AMOC}
As third test for the proposed EWS we consider a higher-dimensional system from ocean circulation modeling.
The Atlantic Meridional Overturnic Circulation (AMOC) is a complex system of currents which transport warm water across the Atlantic ocean, from the Equator to the poles.
In the northern hemisphere, the AMOC forms a closed loop of currents.
The warm and salty surface water of the lower latitudes is transported, by the Gulf Stream and the North Atlantic Current, to the Labrador and Nordic seas.
Once there, the water cools and densifies, sinking from the surface and returning to the Equator as the North Atlantic Deep Water, thus closing the loop.
The strength of these currents is regulated by the balance between evaporation processes at the Equator and ice-melting processes in the North Atlantic. 
A weakened state of the circulation strength of the AMOC can deeply affect the climate on a global scale \cite{Sutton2012,vanWesten025b,Dijkstra2026}.

Modeling such system follows two main approaches in the literature: on the one hand multiphysics General Circulation Models resolve the dynamics of these ocean currents on a global scale and, while useful for climate predictions, they are too complex to analyse \cite{vanWesten2023,vanWesten025a}; on the other hand low-dimensional conceptual models try to capture the essential characteristic of the dynamics of the system by looking at its fundamental mechanics.
While simplistic in nature, conceptual models of the AMOC allow for detailed analysis of the qualitative, long-term behaviour of the system through the lenses of dynamical systems and bifurcation theory \cite{Alkhayuon2019,Mancini2026}.

Among the first of such models, a reduced formulation of the North Atlantic to a 2-box system became seminal in the literature: in this formulation, one box corresponds to the surface water at the Equator while the other correponds to the waters at the North Pole.
These two boxes are connected by fluxes of temperature and salinity which model surface flow as well as deep-water mixing.
\begin{figure}[t]
    \centering 
    \includegraphics[keepaspectratio, width=\textwidth]{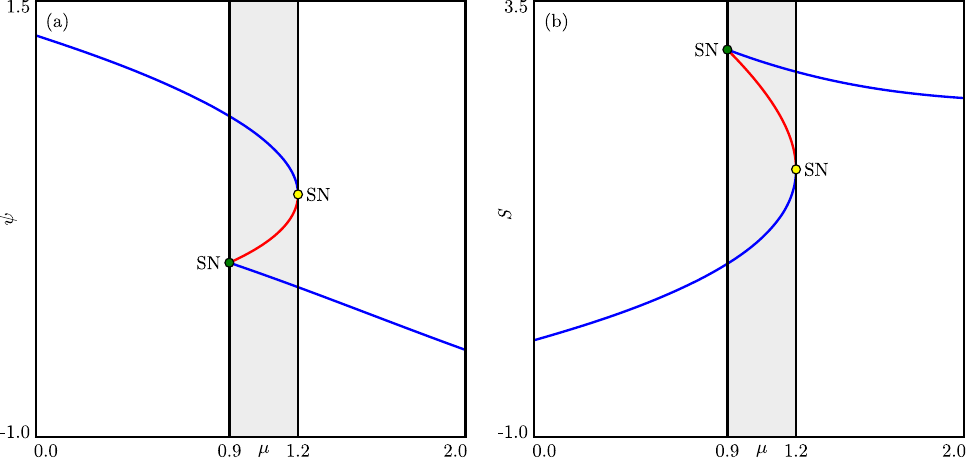}
    \caption{Bifurcation diagrams of the flow-salinity formulation of Stommel's model \eqref{eq:stommel_model_transformed} for $\eta=3$ and $\beta=0.3$.
             Branches of stable equilibria are in blue, while unstable ones are in red.
             The region of bistability (shaded box), between parameter values $\mu=0.9$ and $\mu\approx1.2$, is bounded by a SN bifurcation (yellow circle) and a non-smooth SN bifurcation (green circle).}
    \label{fig:stommel_bif_diag}
\end{figure}
This model, introduced by Stommel in 1961 \cite{Stommel1961}, consists of a two-dimensional dynamical system of the form
\begin{equation}
        \begin{aligned}\label{eq:stommel_model}
                \dot{T} &= \eta - T(1 + |T-S|)\,,\\
                \dot{S} &= \mu - S(\beta + |T-S|)\,,
        \end{aligned}
\end{equation}
where the state variables are the (dimensionless) temperature ($T$) and salinity ($S$) northward gradients while the parameters $\eta$, $\mu$ and $\beta$ model the thermal forcing, freshwater forcing and timescale ratio between thermal and freshwater restoring actions.

Both state variables are proportional to the absolute value of the circulation strength $\psi:=T-S$.
By taking the time derivative of this quantity $\dot{\psi} = \dot{T} - \dot{S}$ and substituting the right-hand sides in Eqs.~\eqref{eq:stommel_model} one can rewrite the two-dimensional system in terms of the circulation strength
\begin{equation}
         \begin{aligned}\label{eq:stommel_model_transformed}
                \dot{\psi} &= \eta - \mu -\psi(1+|\psi|) + S(1-\beta)\,,\\
                \dot{S} &= \mu - S(\beta + |\psi|)\,.
        \end{aligned}
\end{equation}
By fixing $\eta$ and $\beta$ we can study the bifurcation structure of \eqref{eq:stommel_model_transformed} at different regimes of freshwater forcing $\mu$ in the North Atlantic.
Both the circulation strength $\psi$ and salinity $S$ gradients show bistable regions in $\mu$ bounded by two SN bifurcations, one of which is non-smooth \cite{Budd2024} (see the bifurcation diagrams in \fref{fig:stommel_bif_diag}).
In the healthy state, the AMOC's circulation is highly sustained and the advection of heat from the Equator to the North Atlantic is strong.
However, the alternative stable state is one in which the heat transport is weakened by lower values of the circulation strength which causes a collapse in thermal flux and an increased salinity gradient.

The current state of the AMOC is that of high strength, however it has long been established that the increased freshwater forcing of the past two centuries, due mainly to the melting of the Greenland ice sheets, has resulted in the weakening of the heat advection \cite{Dakos2008,Boers2021b,Lohmann2021}.
Given the scale at which a collapse of the AMOC to its weakened state affects the global system, a multitude of scientific communities have produced efforts in the study of the catastrophic tipping event as well as in its prediction \cite{Boers2021a,vanWesten2024,Mehling2024,Chapman2025,Lohmann2025a,Lohmann2025b,Dijkstra2026}.
We contribute to these endeavours by exploiting the interpretability of the modified escape EWS to assess whether the tipping of the AMOC to its weakened state is detectable from multivariate solutions of a high-dimensional system.

\begin{figure}[t]
    \centering 
    \includegraphics[keepaspectratio, width=\textwidth]{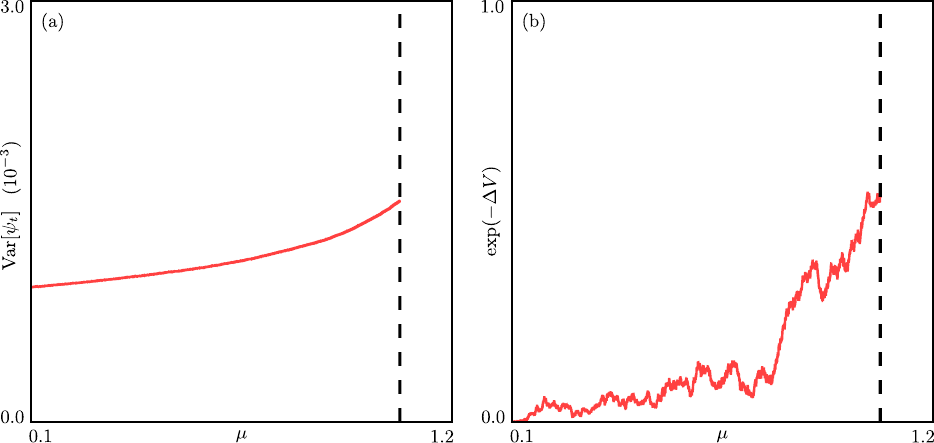}
    \caption{Comparison between the increase in variance and modified escape EWS applied to the one-dimensional circulation observable $\psi$ of the temperature and salinity gradients in the North Atlantic.
            We fix parameters $\eta=3$ and $\beta=0.3$ while the freshwater forcing $\mu$ is varied between $0$ and $1.2$ with rate $\varepsilon=5\cdot10^{-5}$ (see \fref{fig:stommel_bif_diag} for bifurcation diagrams of the stationary system).
            We simulate $100$ sample paths of the slow-fast form of \eqref{eq:stommel_model} and compute the ensemble mean of the modified escape EWS and variance (both in red) across a sliding window of $10\%$ the length of each sample until the earliest tipping (dashed line).
            Notice that \eqref{eq:stommel_model} is not a gradient system (the right-hand side is not a conservative vector field), hence no analytical ground truth of the modified escape EWS is available (panel (b)).
            Furthermore the system is two-dimensional, hence no closed form stationary approximation of the FPE is available, explaining the lack of an analytical reference for the increase in variance (panel (a)).
    We once more remark the different range of values for the variance of $\psi_t$ compared to the previous two cases of vegetation biomass and SN normal form (see Figures \ref{fig:ramped_SN_EWS}(c1) and \ref{fig:ramped_may_ews}(a) respectively).}
    \label{fig:ramped_stommel_ews}
\end{figure}
Two approaches are available to us; we either solve Stommel's $2-$box model in its changed coordinate system \eqref{eq:stommel_model_transformed} and then apply the proposed technique to the circulation strength $\psi$ or; we solve the more physical, original system \eqref{eq:stommel_model} in temperature $T$ and salinity $S$ and then compute the circulation strength $\psi=T-S$ as an \textit{observable} of our multivariate solution.
The two approaches are identical in principle, however we hereby present results based on the latter as it has a more direct application to a real-world scenario where measuring temperature and salinity gradients in the North Atlantic is easier than measuring the strength of its circulation directly.

As we have done for the previous two cases, we formulate a slow-fast SDE of the form \eqref{eq:fast_slow_sde} with $X_t = (T, S)$, the drift term $f(X_t,\mu)$ being the right-hand side of \eqref{eq:stommel_model} and the freshwater forcing $\mu$ evolving slowly with linear ramp $\varepsilon\ll1$.
We compute an ensemble of $100$ sample paths solving the slow-fast system and then, out of the bivariate solutions, we compute a univariate ensemble of timeseries representing the circulation observable $\psi$.
Equivalently to the previous cases, we then compute the ensemble averages of timeseries variance and modified escape EWS, which we report in \fref{fig:ramped_stommel_ews}.

Remarkably, the modified escape EWS reconstructed from the MLE of the cubic potential retains its upward trend close to the bifurcation even though the underlying system \eqref{eq:stommel_model} does not have a confining potential (\fref{fig:ramped_stommel_ews} (b)).
We conclude that the modified escape EWS provides a much more confident indication for the tipping of the AMOC to its weakened state than the traditional increase in variance (\fref{fig:ramped_stommel_ews} (a)), at least for the conceptual, $2$-box model considered in this paper.
\section{Performance analysis}\label{sec:analysis}
The three test cases of \sref{sec:applications} demonstrate that the modified escape EWS can be estimated from slowly varying timeseries following the MLE of the potential discussed in \sref{sec:mle}.
While those examples show the applicability of Algorithm~\ref{alg:modified_ews} to different systems approaching B-tipping, they all refer to similar values of the rate of change parameter $0<\varepsilon\ll1$, noise intensities $0<\sigma<1$, sampling frequencies $0<\Delta T<1$ and timeseries lengths $N\to\infty$.
In the following section we investigate how the performance of the proposed inference method is affected by varying these parameters.
We thus perform statistical analysis of the relevant quantities involved in the estimation of the modified escape EWS following Algorithm~\ref{alg:modified_ews}.
As a representative example, we consider a one-dimensional, slow-fast system of the form
\begin{equation}
        \begin{aligned}\label{eq:bistable_model}
                dX_t &= (-\mu + 4x - 4x^3)dt + \sigma dW_t\,,\\ 
                d\mu &= \varepsilon dt\,,
        \end{aligned}
\end{equation}
where the bistable drift $f(x,\mu)= -\mu + 4x - 4x^3$ undergoes a SN bifurcation at $\mu_c\approx1.53$.
This choice preserves the same characteristics of the three cases considered in the previous section (slow ramp through B-tipping in a regime of bistability) while providing yet another example of application of the proposed methodology.
The deterministic, linear ramp of the bifurcation parameter is written explicitly as the evolution law $\mu(t) = \mu_0 + \varepsilon\,t$ to emphasize the dependence on the starting value $\mu_0 \ll \mu_c$.
Generically, for given sampling frequency $\Delta T$ and timeseries length $N$, the end timestep $T=N\Delta T$ uniquely identifies the parameter value at the end of the ramp $\mu(T)=\mu_0+\varepsilon\,N\Delta T$ for a specific rate $\varepsilon$. 

\subsection{Robustness of the residuals}\label{subsec:residuals}

\begin{figure}[t]
    \centering 
    \includegraphics[keepaspectratio, width=\textwidth]{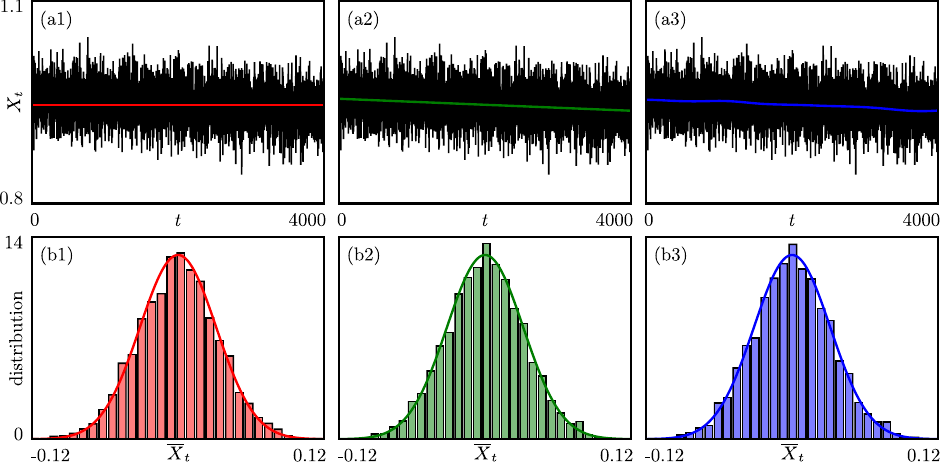}
    \caption{Comparison of detrending techniques for a solution $\mathcal{X}$ of \eqref{eq:bistable_model} with $\mu_0=0$, $\varepsilon=5\cdot10^{-4}$ and $\sigma=0.1$; the simulation stops at the parameter value $\mu(T)=0.2$.
            Panels (a1-a3): trends for the same ramped timeseries (black curves) are computed using the sample mean (red curve), linear fit (green curve) and EMD (blue).
            Panels (b1-b3): comparison between the stationary density of the OUP approximating $\mathcal{X}$ (colored curves) and the empirical distribution of the residuals obtained by sample mean removal (red histogram), linear fit (green histogram) and EMD (blue histogram).}
    \label{fig:residuals_distribution}
\end{figure}

One of the assumptions of the MLE of the modified escape EWS is that the underlying distribution of the sample path is stationary.
However, solutions $\mathcal{X} = \{X_0,\dots,X_N\}$ of slowly ramped systems, such as \eqref{eq:bistable_model}, only satisfy this condition in the singular limit $\varepsilon\to0$ and for long enough time horizons $N\to\infty$.
For strictly positive rates $\varepsilon>0$ the deterministic drift $f$ is nonautonomous and thus the FPE of the system does not admit a Boltzmann-like stationary solution \eqref{eq:invariant_pdf}.
Nonetheless, by considering slow-fast regimes for which $0<\varepsilon\ll1$, we can obtain a quasi-stationary sample $\overline{\mathcal{X}}=\{\overline{X}_1,\dots,\overline{X}_N\}$ by removing the underlying trend in the ramped sample path $\mathcal{X}$.
The detrending of non-stationary timeseries (see step 1. in Algorithm~\ref{alg:modified_ews}) is a fundamental pre-processing stage in the derivation and analysis of EWS as the quality of the extracted quasi-stationary residuals can affect the performace of the estimation.

In the test cases considered throughout \sref{sec:applications} we obtained the residuals by employing the Empirical Mode Decomposition (EMD) \cite{Huang1998} which is a data-driven technique that extracts the underlying trend of a timeseries through a sifting procedure.
Specifically, EMD iteratively decomposes a timeseries in so-called Intrinsic Mode Functions (IMFs), each capturing the different scales of fluctuations in the timeseries.
The trend of a sample path $\mathcal{X}$ thus corresponds to the IMF with the largest scale of fluctuation (see \cite{Huang1998} for details) and the residuals $\overline{\mathcal{X}}$ are obtained by removing the leading IMF from the sample path.

A desirable property of EMD is that it does not require fine calibrations as the IMFs are extractred through the envelopes of the timeseries.
This makes the application of EMD particularly suitable in the computation of quasi-stationary residuals $\overline{\mathcal{X}}$ as it readily extends to univariate timeseries with generic trends without the need of fine-tuning additional parameters.
Similar data-driven techniques, often used in applications, are centering (i.e. removal of the sample mean from the timeseries) and linear fit (i.e. removal of the slope from the timeseries).

We simulate an ensemble of $1{,}000$ sample paths of \eqref{eq:bistable_model} ramped between $\mu_0=0.0$ and $\mu(T)=0.2$ with $\varepsilon=5\cdot10^{-4}$; for each of those sample paths we apply detrending by EMD, linear fit and centering and analyse the residuals separately.
In order to assess the quality of the obtained residuals with respect to the stationary assumption, an OUP, approximating the slowly ramped sample paths, is constructed with coefficients $\alpha=f'(x=a)$ and $a$ being evaluated at $\mu(T)=0.2$ (see Eq.~\eqref{eq:oup}).

\fref{fig:residuals_distribution} shows a representative example of this experiment where a sample path of the simulated ensemble (black curves in panels (a1-a3)) is detrended by centering (red curve, panel (a1)), linear fit (green curve, panel (a2)) and EMD (blue curve, panel (a3)).
The distribution of the quasi-stationary residuals (histograms in panels (b1-b3)) are then compared to the stationary distribution of the OUP (see Eq.~\eqref{eq:oup_density}).
To quantify the deviation of the residuals from the stationary assumption we compute the relative error between the residuals' statistics $\overline{q}_k$ and the central moment $q_k$ of the OUP
\begin{equation}\label{eq:relative_error}
        E_k\coloneqq \frac{|\overline{q}_k - q_k|}{|q_k|}\,,
\end{equation}
where $k=1,2,3$ denote the first (mean), second (variance) and third (skewness) moment of the distribution.
The results of this experiment, compiled in Table~\ref{tab:relative_error}, show that, within the slow ramp regime $0<\varepsilon\ll1$, the three detrending techniques largely preserve the stationary assumption in their residuals with average relative errors well within $10\%$.
Notice that, while the three techniques yield similar residuals, EMD has the advantage of generalising to more complex and nonlinear transient regimes beyond the slow ramp forcing considered throughout this manuscript.
\begin{table}[t]
  \centering
  \caption{Relative errors statistics \eqref{eq:relative_error} of an ensemble
    of $1{,}000$ residuals obtained by different detrending techniques.}
  \label{tab:relative_error}
  \small
  \setlength{\tabcolsep}{3pt}
  \begin{tabular*}{\textwidth}{@{\extracolsep{\fill}} l ccc ccc ccc @{}}
    \toprule
    & \multicolumn{3}{c}{Centering}
    & \multicolumn{3}{c}{Linear fit}
    & \multicolumn{3}{c}{EMD} \\
    \cmidrule(lr){2-4} \cmidrule(lr){5-7} \cmidrule(lr){8-10}
    & min & mean & max & min & mean & max & min & mean & max \\
    \midrule
    $E_1$ & $5\cdot10^{-19}$ & $4\cdot10^{-17}$ & $4\cdot10^{-16}$
          & $5\cdot10^{-19}$ & $1\cdot10^{-16}$ & $6\cdot10^{-16}$
          & $2\cdot10^{-6}$  & $3\cdot10^{-4}$  & $2\cdot10^{-3}$ \\
    $E_2$ & $0.02$ & $0.03$ & $0.05$ & $0.03$ & $0.04$ & $0.05$
          & $0.03$ & $0.04$ & $0.05$ \\
    $E_3$ & $2\cdot10^{-5}$ & $0.03$ & $0.16$ & $6\cdot10^{-5}$ & $0.03$ & $0.19$
          & $1\cdot10^{-4}$ & $0.03$ & $0.18$ \\
    \bottomrule
  \end{tabular*}
\end{table}

\subsection{Accuracy and sensitivity}\label{subsec:sensitivity}
In the singular limit, the MLE of the modified escape EWS converges, in mean, to the ground truth as the underlying system is driven closer to B-tipping (see \fref{fig:stationary_SN_ews}).
Specifically, as $N\to\infty$, the LLS solutions of the potential estimation are normally distributed with their standard deviations decreasing as the bifurcation parameter $\mu$ approaches the SN bifurcation $\mu_c$ (see \fref{fig:lls_solutions}).
The statistical properties of the LLS solutions are prescribed by the normally distributed uncertainties $\Xi\sim\mathcal{N}(0,\sigma^2/\Delta t)$ in the first-order increments (see Eq.~\eqref{eq:increments}).
As a result, the accuracy of the LLS solutions, and hence the MLE of the modified escape EWS, are reflected by the specific choices of additive noise intensity $\sigma$ and sampling frequency $\Delta t$.
Furthermore, it is reasonable to expect that shorter (i.e. smaller $N$) and faster (i.e. larger $\varepsilon$) timeseries would produce significantly different residuals and thus hinder the stationary assumption discussed above. 

\begin{figure}[t]
    \centering 
    \includegraphics[keepaspectratio, width=\textwidth]{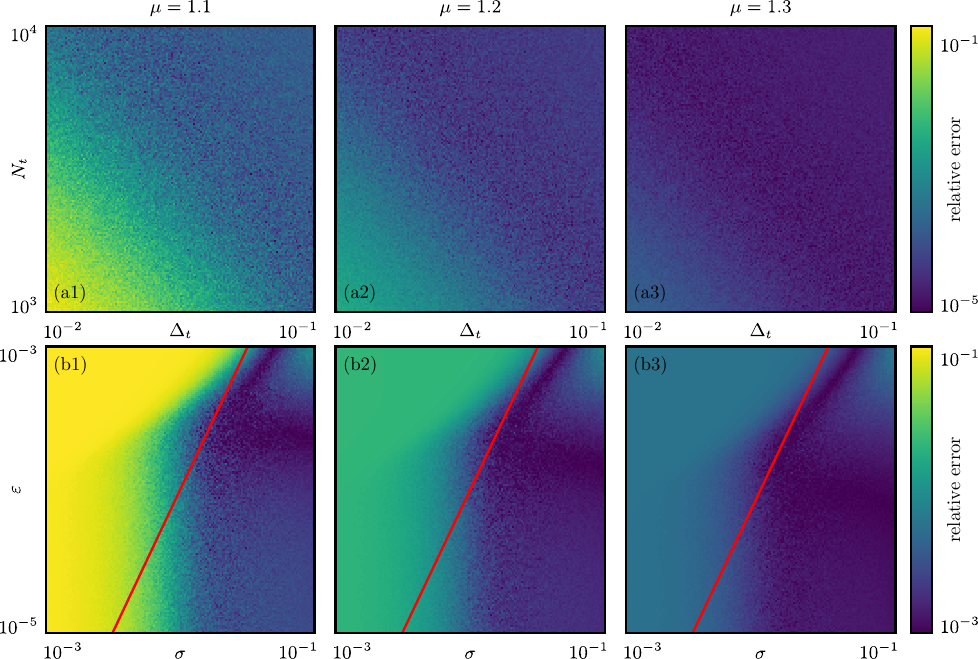}
    \caption{Sensitivity of the relative error of the MLE modified escape EWS under different regimes.
            Each pixel in the heatmap corresponds to the ensemble mean relative error of $1{,}000$ sample paths of \eqref{eq:bistable_model}.
            Panels (a1-a3): in the stationary case ($\varepsilon=0$) and with fixed small noise ($\sigma=0.05$) longer timeseries, sampled at larger frequencies, yield a more accurate estimation of the modified escape EWS.
            Panels (b1-b3): in the slow-fast regime, for fixed timeseries length $N=10{,}000$ and sampling frequency $\Delta t = 0.1$, the accuracy of the estimated modified escape EWS depends, non-trivially, from the level of additive noise in the system with smaller noise seemingly yielding less accurate estimations.
            The red curve represents the $\sigma=\sqrt{\varepsilon}$ threshold.}
    \label{fig:statistical_performance}
\end{figure}

We investigate the sensitivity of the MLE of the modified escape EWS when the above constraints (long and slow timeseries with small noise and large sampling frequency) are relaxed.
Specifically, we test the accuracy of Algorithm~\ref{alg:modified_ews} in two different settings:
\begin{itemize}
     \item by varying the timeseries length $N$ and sampling frequency $\Delta t$ in the singular limit $\varepsilon\to0$ and in a regime of small noise $\sigma\to0$ and;
     \item by varying the parameter's rate of change $\varepsilon$ and noise level $\sigma$ with fixed timeseries length $N$ and sampling frequency $\Delta t$.
\end{itemize}
In both cases we run an ensemble simulation of $1{,}000$ sample paths of \eqref{eq:bistable_model} and evaluate the accuracy of the proposed method by computing the a-posteriori relative error (see Eq.~\eqref{eq:relative_error}) between the estimated modified escape EWS and the ground truth.
In order to assess whether the convergence in mean persists in different regimes, we compare the relative errors at different parameter values as the system is driven closer to B-tipping (i.e. as $\mu\to\mu_c$).
 
The results of this experiment are depicted in \fref{fig:statistical_performance}.
In the stationary case (panels (a1-a3)) we observe that the modified escape EWS estimate is more accurate for longer timeseries $N\to10^4$ and larger sampling frequencies $\Delta t\to10^{-1}$.
This is not surprising considering that the uncertainty in the LLS regression is inversely proportional to $\Delta t$.
Furthermore, we notice that the increase in accuracy close to B-tipping is preserved for varying $N$ and $\Delta t$ as the mean relative error for the system far from the SN bifurcation $\mu=1.1$ (panel (a1)) is overall larger than the one for the system closer to the SN bifurcation $\mu=1.3$ (panel (a3)).
In the slow-fast regime (panels (b1-b3)) we keep the timeseries length $N=10^4$ and sampling frequency $\Delta t=10^{-1}$ fixed while we generate solutions of \eqref{eq:bistable_model} at varying rates $\varepsilon$ and noise intensities $\sigma$.
For each $\varepsilon$ we keep $\mu(T)$ fixed and adjust $\mu_0$ so that each sample path contains exactly $N=10^4$ realizations.
In particular $\mu(T)\in\{1.1,1.2,1.3\}$ in panels (b1), (b2) and (b3) respectively; these parameter values are then used to compute the ground truth of the modified escape EWS.
Intuitively, we observe that faster rates $\varepsilon\to10^{-3}$ yield overall less accurate estimations due to the residuals deviating further from the stationary assumption (as discussed in the previous section).
On the other hand, the estimations become less accurate for smaller noise levels $\sigma\to10^{-3}$; while Kramer's law holds in regimes of small noise, the LLS estimation of the confining potential requires larger excursions around the stable equilibrium to better infer the local shape of the basin of attraction.
In particular, we notice that the non-trivial decay of the relative error shows a complex relationship modelling the interplay of $\varepsilon$ and $\sigma$ in different regimes, meaning that slower and less noisy timeseries do not necessarly yield better estimates.
This emphasizes the fact that the MLE of the modified escape EWS is sensitive to perturbations in forced systems and in particular how the threshold $\sigma=\sqrt{\varepsilon}$ (red curves in panels (b1-b3)), separating regimes of slow passage through B-tipping from N-tipping dominated events, provides an actionable guideline for applications.

\subsection{Statistical significance}\label{subsec:significance}
\begin{figure}[!t]
    \centering 
    \includegraphics[keepaspectratio, width=\textwidth]{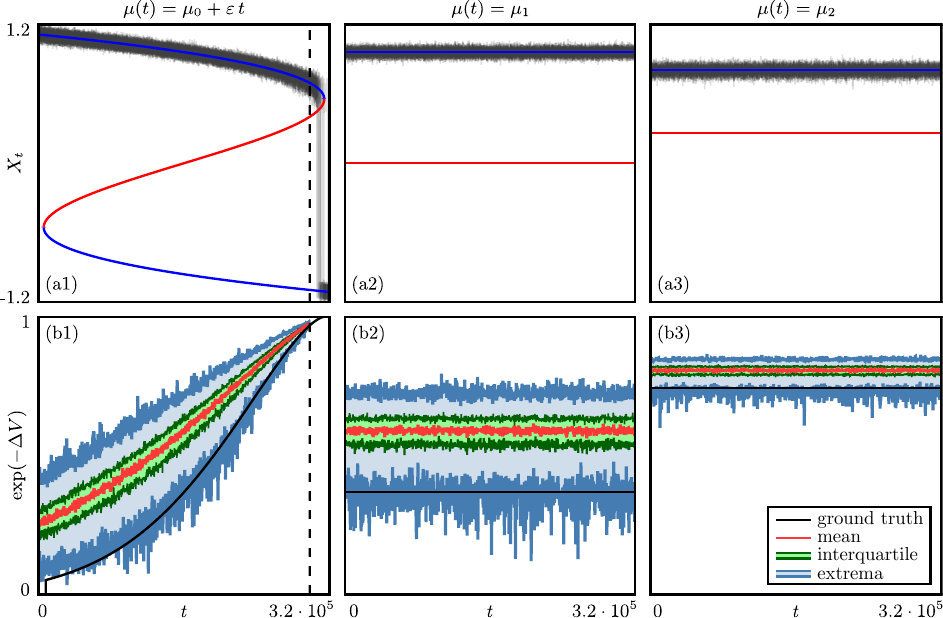}
    \caption{Time-varying estimations of the modified escape EWS for an enseble $1{,}000$ slowly ramped (panels (a1)-(b1)) and stationary (panels (a2-a3)-(b2-b3)) solutions of \eqref{eq:bistable_model}.
            In the ramped case the parameter varies between $\mu_0=-1.6$ to $\mu(T)=1.6$ with rate $\varepsilon=1\cdot10^{-5}$. 
            The sample paths, shown by the black curves in panel (a1), follow the upper stable branch of the bifurcation diagram (blue curve) of the frozen system before tipping across the unstable branch (red curve).
            The earliest tipping in the ensemble is indicated by the vertical dashed line.
            In the stationary cases the system is in the singular limit $\varepsilon\to0$ and thus the solutions, shown by the black curves in panels (a2-a3), are confined around the stable equilibrium (blue line) and do not tip across the unstable equilibrium (red curve) within the simulated time horizon.
            The stationary cases are simulated at parameter values $\mu=0$ (panel (a2)) and $\mu=1$ (panel (a3)).
            In all cases, a rolling window of $10\%$ the length of the sample paths, is used to estimate the modified escape EWS of the entire ensemble (panels (b1-b3)).
            The ensemble envelope (i.e. maximum and minimum values) of the estimations are depicted by the blue shaded regions, whereas the interquartile range is depicted by the green shaded regions. 
            The ensemble mean is shown by the red curves and it is compared with the ground truth of the frozen system, depicted by the black curves.}
    \label{fig:statistical_significance}
\end{figure}
The convergence in mean of the estimated modified escape EWS guarantees accurate warnings when an ensemble of sample paths is available.
In real-world applications however, this is rarely the case as access to i.i.d. samples of a given system is limited by physical and technical constraints.
We therefore investigate how significant the upward trend for the estimated modified escape EWS is for individual realizations of \eqref{eq:bistable_model}.
In particular, we are interested in establishing whether the estimations from single sample paths can lead to missed allarms (i.e. false negatives) and/or spurious warnings (i.e. false positives).

To do so we perform an ensemble simulation of $1{,}000$ sample paths ramped from $\mu_0=-1.6$ to $\mu(T)=\mu_c$ with $\varepsilon=1\cdot10^{-5}$.
We then use Algorithm~\ref{alg:modified_ews} to estimate the modified escape EWS across a sliding window containing $N=10^4$ realizations of each sample path.
We repeat the same procedure for stationary simulations of \eqref{eq:bistable_model} at parameter values $\mu_1=0$ (i.e. at a safe distance from B-tipping) and $\mu_2=1$ (closer to B-tipping but still at a safe distance).
We report the results of this analysis in \fref{fig:statistical_significance}.
The sample paths (black curves in panels (a1-a3)) are computed across a time horizon $T=3.2\cdot10^5$. 
For the slowly ramped case (panel (a1)), the timeseries will all tip, from the drifting stable equilibrium (blue curves), across the unstable threshold (red curve), within the prescribed time horizon; this provides a sample of positive tipping events.
For the stationary cases (panels (a2-a3)) the solutions will stay bounded around a small region of the stable equilibrium (blue line) and they will not tip across the unstable threshold (red line) within the time horizon (red line); this provides a sample of negative tipping events.
In the ramped case, the estimated ensemble mean of the modified escape EWS (red curve, panel (b1)) overestimates the proximity to B-tipping given by the ground truth (black curve) partly due to less amount of data available; notice that the familiar convergence in mean, as $\mu\to\mu_c$, implies that the overestimation bias vanishes close to B-tipping.
Furthermore, while the interquartile range (green shaded region) shows a reasonable bounded region for the estimated signal, the outliers in the ensemble (blue shaded region) deviate significantly from the mean.
This feature of the MLE of the modified escape EWS is also observed in the stationary cases with the deviations reducing as the system is driven closer to B-tipping.

We conclude that, the overestimation bias in the ensemble mean, and the large deviations of the outliers from such mean can lead to spurious warnings when a single trajectory is considered.
However, considering that the interquartile range provides a considerably smaller region around the mean, and that such deviations vanish as $\mu\to\mu_c$, we can also conclude that the overall trends of the estimations become more significant as the system approaches B-tipping, which minimizes the occurence of missed alarms from single trajectories.
Notice that, variability in the outliers of the estimated EWS can be alleviated by considering moving averages of the timeseries, thus yielding much more reliable indicators that can also minimize spurious warnings.

\section{Discussion}\label{sec:discussion}
In this manuscript we have introduced a new EWS for critical transitions in complex systems.
This EWS is a modified version of Kramer's formula, used in statistical mechanics to quantify the rate of escapes of particles across potential barriers.
The main aim of this EWS is to be universally interpretable for the estimation of proximity of B-tipping in natural and human systems. 
\begin{figure}[t]
    \centering 
    \includegraphics[keepaspectratio, width=\textwidth]{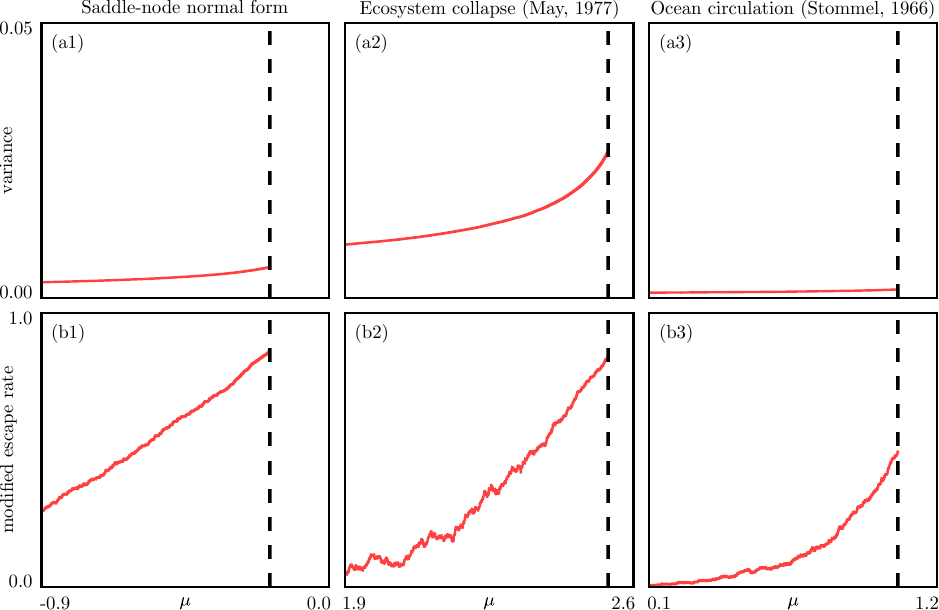}
    \caption{Comparison of the different range in values of the timeseries variance (panels (a1-a3)) and the modified escape EWS (panels (b1-b3)) for the three different test cases described in \sref{sec:applications}: a one-dimensional SN normal form \eqref{eq:slow_fast_SN} (panels (a1,b1)); a one-dimensional model of ecosystem collapse \eqref{eq:may_model} (panels (a2,b2)) and; a one-dimensional observable of a two-dimensional flow-salinity model of ocean circulation in the North Atlantic \eqref{eq:stommel_model_transformed} (panels (a3,b3)).
            The curves in each box represent the ensemble mean of the respective quantities computed for $1,000$ sample paths. 
            The dashed vertical line indicate the earliest tipping in the ensemble of sample paths. 
    Notice how the modified escape EWS yields a comparable range of values for the three different systems, contrary to the increase in variance.}
    \label{fig:ews_comparison}
\end{figure}
This extends the work of previous, case-specific interpretable EWS such as the $2-$sigma criterion \cite{Drake2010} and Kendall-$\tau$ trend \cite{Harris2020} of statistical indicators, the classifier score calibrated on training data \cite{Brett2020} and the generalisation of autocorrelation estimators \cite{Rodal2022}.

The proposed measure has the advantageous property of increasing monotonically in the bounded range $(0,1]$: values close to $0$ indicate strong potential confinement (safe distance from B-tipping); values close to $1$ indicate weak potential confinement (close proximity to B-tipping).
Furthermore, the modified escape EWS introduced in this manuscript is system-independent, meaning that the monotonic increase from $0$ to $1$ is registered for any low-dimensional complex system provided that it has a gradient structure and that its stochastic fluctuations are state-independent (i.e. additive noise).
While Kramer's formula is still valid for regimes of state-dependent fluctuations (i.e. multiplicative noise), the exponential decay no longer depends exclusively from the potential barrier, which further complicates the extraction of a plausible and interpretable EWS.
We remark that the modified escape EWS shows an earlier upward trend and consistent range of values for different systems (\fref{fig:ews_comparison} (b1-b3)), contrary to the canonical increase in timeseries variance which has different ranges for different systems (\fref{fig:ews_comparison} (a1-a3)).
Due to these properties we argue that the proposed modified escape EWS is a more reliable precursor of B-tipping than the traditional diagnostics of increasing sample variance.
 
Additionally, in the present paper we derived a MLE method to extract the modified escape EWS from timeseries without the need of prior knowledge of the underlying mechanism generating it.
We showcased how a simple LLS regression from the given sample can yield accurate enough approximations of the modified escape EWS given its convergence in mean in the limit of large number of samples.
We observed that the approximation becomes more accurate as the system is driven closer to B-tipping, thereby providing a consistent prevention against false negatives.

Extensions to multivariate signals is also a subject of high interest among the various communities trying to unlock EWS of high-dimensional systems using spatio-temporal information \cite{Donovan2015,Rodal2022,Lohmann2024,Dijkstra2026}.
The emphasis of these efforts is to generalize precursors of B-tipping events (e.g. CSD) to spatially-extended models.
        Aggregations of the spatial information into mean-field approximations is likely to inhibit important mechanisms that drive the tipping.
        One proposed methodology has been the choice of optimal, one-dimensioal observables of these spatially-extended systems that still retains sufficient information to capture the precursors \cite{Lohmann2025b}.
        It is not clear at this stage how reliable traditional temporal, spatial and spatio-temporal indicators should perform in generic high-dimensional systems \cite{Robinson2025}.
        In this paper we obtain evidence that applications of the proposed methodology to one-dimensional reductions of higher-dimensional systems can yield confident results, provided an adequate choice of the monitored observable.

Further investigation in improving the accuracy of the technique is required however, as the reconstructed signal from empirical samples seem to be subject to higher uncertanties than the statistical moments counterparts.
While this falls outside the scope of the current manuscript, which instead shows a proof-of-concept of the methodology, as well as introducing a novel and interpretable EWS to the literature, we speculate that refinements of the MLE framework for this problem could enhance its accuracy.
\appendix
\section{Derivation of the modified escape rate as an EWS}\label{appendix_a}
In the following we explain the motivations behind the choices of the definition of a modified escape EWS
\begin{equation*}
        \exp(-\Delta V)\,,
\end{equation*}
as given in \eqref{eq:escape_ews}, rather than using Kramer's original formula for the escape rate \eqref{eq:kramer_escape_formula}
\begin{equation*}
        \frac{\sqrt{|V''(b)|V''(a)}}{2\pi}\exp\bigg(-\frac{\Delta V}{D}\bigg)\,.
\end{equation*}
We remark that the modified escape EWS we propose is obtained by Kramer's formula by dropping the prefactor $\sqrt{|V''(b)|V''(a)}/2\pi$ and by setting $D=1$.
To motivate this choice we emphasize that one of the main contributions of our porposed EWS is interpretability.

As explained in the introduction, the increase in variance is only valuable as an EWS of B-tipping if: either one has access to a large enough sample so that monitoring relative changes is possible or; one knows the baseline value of the variance, for the same system, when at a safe distance from the critical transition.
Regardless of the two cases above, the overarching strategy for employing the increase in variance as an EWS is that one needs to compare the detected value from a given sample to another value, of the same system, from a different sample.
A single value of the variance does not provide any insight for the detection of CSD, which is at the core of EWS of B-tipping.

On the contrary, the modified escape EWS is bounded in $(0,1]$, where values close to $0$ immediately translate to safe distance of the system from the critical transition, while values close to $1$ indicate a dangerous proximity to a critical transition.
This advantage of interpretability of the modified escape EWS is further emphasized by the fact that, as a gradient system approaches B-tipping, the potential barrier $\Delta V\to0$ monotonically with $\mu\to\mu_c$.
The modified escape EWS is a monotonic function of $\Delta V$ which means that, at least analytically, if the system is driven closer to B-tipping the modified escape EWS can only increase and viceversa.
This leaves no ambiguity in interpretation for the measured trend of the modified escape EWS. 
The same cannot be said about the increase in variance which, as shown in Figures \ref{fig:stationary_SN_ews} and \ref{fig:ramped_may_ews}, is non-monotonic with $\mu\to\mu_c$.
We conclude that the two main advantages of using the modified escape EWS, as a precursor of CSD, over the increase in variance, are:
\begin{itemize}
     \item interpretability provided by the boundedness of its range in $(0,1]$;
     \item monotonicity of increase with respect to the system's \textit{``distance''} from B-tipping.
\end{itemize}
\begin{figure}[t]
    \centering 
    \includegraphics[keepaspectratio, width=\textwidth]{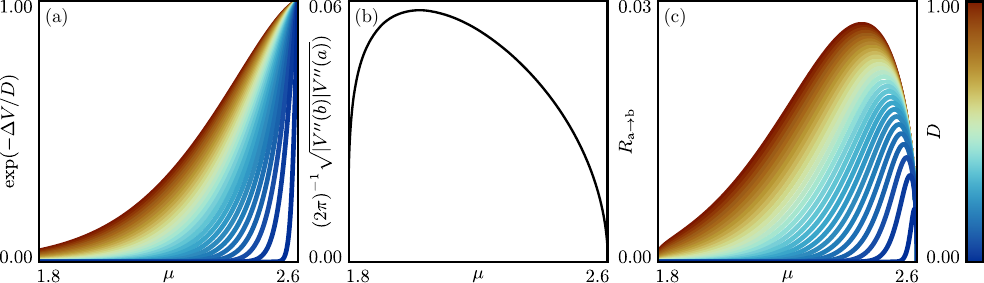}
    \caption{Comparison of the modified escape EWS (panel (a)) and Kramer's escape rate \eqref{eq:kramer_escape_formula} (panel (c)) at different values of the diffusion level $D$ (color-coded in the colorbar on the right).
            The prefactor term of Kramer's escape rate (panel (b)) is independent on $D$ and thus shown by a single black curve.}
    \label{fig:kramers_formula}
\end{figure}

Due to the above, we show that including the dependence of the proposed EWS on the stochastic diffusion term $D$ and the curvature-depending prefactor $\sqrt{|V''(b)|V''(a)}/2\pi$ undermines both the advantages outlined above from an analytical and computational perspective.

First notice that while $\Delta V = V(b) - V(a)\to0$ monotonically with $\mu\to\mu_c$ this is not necessarly true for $(2\pi)^{-1}\sqrt{|V''(b)|V''(a)}$ as shown in \fref{fig:kramers_formula} (b).
While the curvature of both the well $V''(a)$ and hilltop $V''(b)$ of the potential approach $0$ monotonically, as a SN bifurcation involving $x=a$ and $x=b$ is approached, their product does not necessarly.
Furthermore, while it is clear that $\exp(-\Delta V/D)$ is bounded in $(0,1]$, the same cannot be said for the prefactor.
Notice in fact that $\sqrt{|V''(b)|V''(a)}/2\pi$ has range $[0,+\infty)$ for any value of $\mu$ that precedes the SN bifurcation.
Including the prefactor term in our modified escape EWS will thus immediately erase both advantages described above (see \fref{fig:kramers_formula} (c)): it will no longer be bounded in $(0,1]$ but rather in $[0,+\infty)$ (just like the variance) and it will also increase non-monotonically as the SN critical transition is approached (again, just like the variance).

Second, notice that the stochastic diffusion $D$ does not appear in the prefactor and it is only an argument of the exponential part of Kramer's formula.
In particular since $D\in(0,+\infty)$ the exponential $\exp(-\Delta V/D)$ becomes flatter for values of $\mu$ away from B-tippig (see \fref{fig:kramers_formula} (a)) as $D\to0$.
This is undesirable as an EWS because it hinders the \textit{``earliness''} interpretation of the increase of the quantity as $\mu\to\mu_c$.
A good EWS should not only be monotonic in the approach of B-tipping but it should also be \textit{``early''} enough so that preventive actions can be taken to avoid the critical transitions.

This analytical argument in dropping the dependence of the modified escape EWS on $D$ is further corroborated by a simple computational argument.
In the MLE we proposed in \sref{subsec:euler} we estimate the polynomial coefficients $\theta^* = (\theta_1^*,\theta_2^*,\theta_3^*)$ of a cubic potential that better describes the sample $\mathcal{X}=\{X_0,\dots,X_N\}$.
Because of the linearity of the problem on the polynomial coefficients $\theta^*$, the MLE reduces to a LLS problem in $\mathbb{R}^{3}$.
If we were to estimate $D$ as well (which is, in principle, unknown), then the search space becomes $\mathbb{R}^{4}$ and the residual will no longer be linear in the parameter vector (see Eqs. \eqref{eq:increments} and \eqref{eq:residual}).

This increase in complexity for the MLE does not justify the inclusion of the missing information $D$ from the modified escape EWS since it does not change the advantages of boundeness and monotonicity discussed above.
As a result we set $D=1$ and obtain our modified escape EWS to be \eqref{eq:escape_ews}.
\section{Variance-bias effects of Tikhonov regularization of the LLS solutions}\label{appendix_b}
\begin{figure}[!tp]
    \centering 
    \includegraphics[keepaspectratio, width=\textwidth]{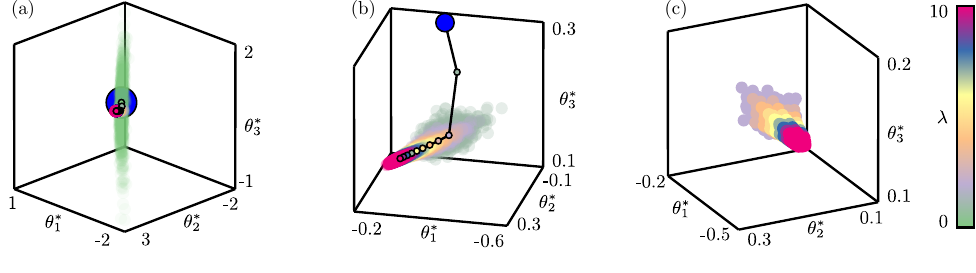}
    \caption{Demonstration of the variance-bias effect of the penalizer $\lambda||\theta||_2$ for Tikhonov regularization of the LLS solutions \eqref{eq:rlls_solutions} at different values of the penalizing parameter $\lambda\in[0,10]$.
            For each value of $\lambda$ we simulate an ensemble of $1,000$ sample paths solving a stationary SN normal form \eqref{eq:slow_fast_SN} at parameter value $\mu=-0.3$ and with $\varepsilon=0$.
            For each sample path we compute the solutions $\theta^*_\text{reg}$ of the regularized LLS problem \eqref{eq:rlls_solutions} and plot its location in the search space $\mathbb{R}^{3}$ with a circle.
            In panel (a) we show the result of this experiment in a cubic subset of $\mathbb{R}^{3}$; the solutions $\theta^*_{\text{reg}} = \theta^*_{\text{reg}}(\lambda)$ are color-coded by the value of the penalizing parameter $\lambda$ used in \eqref{eq:rlls_solutions} (see colorbar on the right).
            The ensemble mean $\hat{\theta}^*_{\text{reg}}(\lambda)$ of the solutions $\theta^*_{\text{reg}}(\lambda)$ are highlighted by a black outline.
            The larger, blue circle represents the location of the ground truth point $\theta = (\mu,0,1/3)$ in the search space.
            In panel (b) we show a zoomed-in detail of panel (a) focused on the solutions for the middle-to-large values of the penalizer $\lambda$.
            In panel (c) we zoom-in further to focus on the solutions for the larger values of $\lambda$.
            The three plots show how, for large values of the penalizer $\lambda$, the variance of the solutions reduces at the cost of introducing undesirable biases in the mean of the solutions.
    Another way of interpreting this result is that the regularized LLS solution \eqref{eq:rlls_solutions} converges in mean to a biased estimator $\tilde{\theta}:=\theta+\delta\theta$ in the limit of large penalizers, i.e. $\hat{\theta}^*_{\text{reg}}(\lambda)\to\tilde{\theta}$ as $\lambda\to\tilde{\lambda}\gg0$.}
    \label{fig:tikhonov_regularization}
\end{figure}
Regularization techniques such as Tikhonov regularization \cite{Tikhonov1963} are employed in ill-conditioned problems often with the aim of reducing the uncertainty in the solutions \cite{Cheng2011}.
They achieve this by introducing a penalizing term to the minimization problem \eqref{eq:quadratic_minimization_problem} 
\begin{equation}\label{eq:rlls_problem}
        \theta^*_{\text{reg}} = \argmin_{\theta\in\mathbb{R}^3}\big(||A\theta - y||_2 + \lambda||\theta||_2\big)\,,
\end{equation}
where $\lambda\geq0$ is a parameter that calibrates the degree of the penalization.
The solution of \eqref{eq:rlls_problem} is
\begin{equation}\label{eq:rlls_solutions}
        \theta^*_{\text{reg}} = (A^{T}A + \alpha\text{ Id}_N)^{-1}A^{T}y\,,
\end{equation}
which, differently from the normal equation \eqref{eq:lls_solution} of the ordinary LLS solution, it adds the penalizing parameter $\lambda$ to the diagonal entries of the square matrix $A^TA$.
Notice that the case $\lambda=0$ reduces to an ordinary LLS problem \eqref{eq:lls_problem}-\eqref{eq:lls_solution}.

The presence of this additional term acts as a penalizer of large solutions $\theta^*\in \mathbb{R}^{3}$ in the $2-$norm sense.
Intuitively, strictly positive values of $\lambda$ reduce the variance of the solutions \eqref{eq:rlls_solutions}, thus yielding more accurate estimations of the polynomial coefficients of the cubic potential.
However, it is well-established that positive values of the penalizer also introduce biases in the estimator \cite{Hoerl1970}.

Specifically, a large value of $\lambda>0$ will succesfully reduce the uncertainty of the solutions \eqref{eq:rlls_solutions} but at the cost of shifting the mean value of those solutions to a spurious point in the search space.
This effect is also known as the variance-bias tradeoff in machine learning. 
In \fref{fig:tikhonov_regularization} we demonstrate such effect by solving regularized LLS estimations \eqref{eq:rlls_solutions} for an ensemble of $1,000$ sample paths of a SN normal form \eqref{eq:slow_fast_SN} in the stationary regime ($\varepsilon=0$) and for increasing values of $\lambda\in[0,10]$.

\printbibliography

\end{document}